\RequirePackage{ifdraft}
\documentclass[acmtog]{acmart}
\acmSubmissionID{750}

\usepackage{booktabs}
\usepackage{threeparttable}
\usepackage{enumitem}
\usepackage{siunitx}
\usepackage[ruled]{algorithm2e}
\usepackage{wrapfig}

\usepackage{amsmath}
\DeclareMathOperator*{\argmin}{arg\,min}
\usepackage[mathscr]{euscript}
 \let\mathscr\relax% just so we can load this and rsfs
\usepackage{amsthm}
\usepackage{mathrsfs}
\usepackage{multicol}
\usepackage{comment}
\usepackage{breqn}
\usepackage{mathtools}
\usepackage{overpic}
\mathtoolsset{showonlyrefs=true}
\newtheorem{remark}{Remark}

\newtheorem{theorem}{Theorem}
\newtheorem{proposition}{Proposition}

\graphicspath{{figures/}}

\def\n{\boldsymbol{n}}

\def\a{\alpha}
\def\b{\beta}
\def\c{\theta}
\def\d{\tau}
\def\e{\mathbf{e}}

\def\eps{{\epsilon}}

\def\toe{{\tilde{\Omega}^\epsilon}}
\def\u{\boldsymbol{u}}
\def\v{\boldsymbol{v}}
\def\A{\boldsymbol A}

\def\H{\boldsymbol H}
\def\R{{\mathbb R}}

\def\T{\mathbb T}

\def\t#1{{\tilde{#1}}}

\def\tw{{\t{\omega}}}

\def\itd{{\,\mathrm{d}}}
\def\fa{\text{ for all }}

\def\bs#1{\boldsymbol{#1}}

\def\pr#1{\left(#1\right)}
\def\veps{\varepsilon}

\def\br#1{\left\{#1\right\}}
\def\pr#1{\left(#1\right)}

\def\tr{\text{tr\,}}
\def\inn#1{\langle #1 \rangle}
\def\grad{\text{grad\,}}

\newcommand{\myHearts}[2]{\color{#1}{\heartsuit}\kern-2.5pt\color{#2}{\heartsuit}}
\usepackage{accents}
\newcommand*{\dt}[1]{%
	\accentset{\mbox{\large\bfseries .}}{#1}}

\acmJournal{TOG}

\copyrightyear{2025}
\acmYear{2025}
\begin{document}

\title{Asymptotic analysis and design of linear elastic shell lattice metamaterials}

\author{Di Zhang}
\email{zhd9702@mail.ustc.edu.cn}
\orcid{0000-0003-0176-3424}
\affiliation{%
  \institution{University of Science and Technology of China}
  \city{Hefei}
  \country{China}
}

\author{Ligang Liu}
\authornote{Corresponding author.}
\affiliation{%
  \institution{University of Science and Technology of China}
  \city{Hefei}
  \country{China}
  }
\email{lgliu@ustc.edu.cn}

\renewcommand{\shortauthors}{Zhang et al.}

\begin{abstract}

We present an asymptotic analysis of shell lattice metamaterials based on Ciarlet’s shell theory, introducing a new metric—\emph{asymptotic directional stiffness} (ADS)—to quantify how the geometry of the middle surface governs the effective stiffness. We prove a convergence theorem that rigorously characterizes ADS and establishes its upper bound, along with necessary and sufficient condition for achieving it. As a key result, our theory provides the first rigorous explanation for the high bulk modulus observed in Triply Periodic Minimal Surfaces (TPMS)-based shell lattices. To optimize ADS on general periodic surfaces, we propose a triangular-mesh-based discretization and shape optimization framework. Numerical experiments validate the theoretical findings and demonstrate the effectiveness of the optimization under various design objectives. 
Our implementation is available at \url{https://github.com/lavenklau/minisurf}.

\end{abstract}




\setcopyright{acmlicensed}
\acmJournal{TOG}
\acmYear{2025} \acmVolume{44} \acmNumber{4} \acmArticle{}
\acmMonth{8}\acmDOI{10.1145/3730888}

\maketitle

\section{Introduction}
\label{sec:introduction}

The application of Triply Periodic Minimal Surfaces (TPMS) in metamaterial design has achieved great success in recent decades~\cite{tpms-survey}.
Numerous simulation and experimental studies have demonstrated their superior stiffness and strength properties~\cite{pnas-shell-limit,D0MH01408G,CHEN2019108109}.
However, these findings are primarily empirical, derived from studies on a limited number of TPMS configurations.
The fundamental mechanisms underlying their exceptional performance are still not fully understood, and it remains unclear whether these advantages generalize to all TPMS configurations.

Since TPMS typically serve as the middle surface of the shell lattice metamaterials~\cite{pnas-shell-limit,opencell2023,MA2022110426, tpms-survey},
understanding how middle surface geometry influences effective stiffness provides a natural starting point for  theoretical investigation.
To this end, we introduce a novel metric called the \emph{asymptotic directional stiffness} (ADS),
defined as the limiting ratio of stiffness to volume fraction as shell thickness approaches zero.
This metric captures the intrinsic relationship between surface geometry and effective stiffness, independent of thickness effects.
By extending  key results from Ciarlet's shell theory~\cite{CiarletGenmem1996}, we derive a convergence theorem that rigorously formulates ADS.
Based on this theorem, we establish an upper bound for ADS, as well as  the necessary and sufficient condition to achieve this bound. 
This condition implies the optimal asymptotic bulk modulus (i.e., the ADS under hydrostatic strain) of TPMS  as a special case,  thus explains previous experimental discoveries~\cite{D0MH01408G}.

For practical applications, we develop a numerical approach to enhance the ADS of middle surfaces through shape optimization.
Our method features a specialized triangular mesh discretization scheme that handles diverse surface topologies and shapes while supporting various optimization objectives.
Unlike traditional shell structure optimization methods~\cite{MA2022110426,Shimoda2014-ao,Ramm01011993}, our formulation decouples the geometry of the middle surface from thickness, enabling the direct search for optimal surface shapes without prescribing thickness-related parameters.
Our  contributions are summarized as follows
\begin{itemize}
	\item  We introduce the asymptotic directional stiffness (ADS) for shell lattice metamaterial and establish the convergence theorem that formalizes it.
	\item We establish the upper bound of ADS and the necessary and sufficient condition to achieve the bound. This result is used to justify the optimal bulk modulus of TPMS.
	\item A novel discretization scheme on triangular mesh that allows us to evaluate and optimize ADS.
\end{itemize}

\section{Related work}

\paragraph{Linear elastic mechanical metamaterial}
The foundation for the design of linear elastic metamaterials is the homogenization theory, which determines the effective properties of the metamaterial by solving a  \emph{cell problem} on the representative volume element (RVE).
Notable  monographs on this topic include~\cite{Allaire2002,Milton_2002}. 
In the graphics community, early works aim to construct microstructure families with prescribed elastic properties by interpolation from samplings~\cite{meta-contr-el2015} or through optimization on a set of parameterized geometry~\cite{elastic-texture,worst-case-stress}.
The constructed microstructure family achieves a wide range of elastic properties.
Besides the regular structures,  other approaches explore the irregular foam structures, either based on Voronoi diagram~\cite{Martinez2016,Martinez2018} or generated from stochastic process~\cite{Martinez2017}.
Following a similar spirit,  \cite{Mart-star-mat} employs the generalized Voronoi diagrams induced by novel star-shaped metrics to improve the aesthetic.
\cite{Schumacher2018} explored the structured sheet materials based on isohedral tilings that covers a wide range of material properties.
Given that these designed structures typically cannot conform to an arbitrary boundary of the model,~\cite{Tozoni2020} introduced the rhombic microstructure with far greater flexibility in shape approximation.
To enlarge the design space, \cite{Makatura-proce} developed a procedural framework to model microstructures by combining various constitutive elements, such as straight/curved beams, thin shells, and solid bulk materials.
A similar approach is seen in~\cite{Liu2024-multilayer}, which integrates truss, plate, and shell structures into a multilayer design strategy.  

\paragraph{Application of TPMS}
The application studies based on TPMS have been exceptionally active in the last decade, we kindly refer the reader to the survey article~\cite{tpms-survey} for a comprehensive overview.
It has become a consensus in the community that TPMS exhibit superior stiffness.
However, such understanding remains largely empirical, originating from numerical simulations and physical experiments on limited types of TPMS~\cite{FENG2021110050,pnas-shell-limit,CHEN2019108109,QURESHI2021121001,Forced-convec,Flow-and-thermal,Torquato2004}.
Due to the lack of theoretical analysis, the underlying reason behind these discoveries is still unclear, and one cannot guarantee that all TPMS, including those yet to be undiscovered, have the same advantageous  properties.
In the graphics community, TPMS have been used to design inner supporting structures in 3D models to enhance stiffness~\cite{hjb-tpms-mech,8703138,10105512} and improve heat dissipation~\cite{hjb-tpms-heat}.
These works either focus on designing structures at a macroscopic scale rather than developing metamaterials or concentrate on limited types of TPMS,  resulting in a lack of general conclusions. 

\paragraph{Asymptotic analysis}
Our analysis of shells is inspired by the serial works~\cite{CiarMem1996, CiarletFlex1996, CiarletGenmem1996}, in which limit equations for three types of shells—membrane shells, flexural shells, and generalized membrane shells—are derived. However, from a rigorous standpoint, directly applying these equations to our problem is not feasible for several reasons.
First, the boundary conditions differ from theirs.
Second, these shell equations rely on specific geometry assumptions on the middle surface that may not be satisfied by an arbitrary periodic surface.
Third, a general periodic middle surface,  due to its non-trivial topology, cannot be covered by a single coordinate chart, which is a common assumption behind these theories.
Thus, technical adaptions of their theoretic results are necessary before the application.
Another crucial distinction in our problem is that we are studying the convergence of a scalar quantity instead of the vector field.
This allows us to establish a very general and concise conclusion without the need to classify different types of shells, as we will show in Section~\ref{sec:asym-ana}.

\paragraph{Simulation and optimization of shell structures}
The existing works can be broadly classified into two categories: parametric and non-parametric methods.
Parametric approaches typically use NURBS~\cite{KIENDL2014148,Kang2016-xx,Hirschler2019-id,CAI2023116218,JIANG2023613} or B\'ezier~\cite{Ramm01011993,Bletzinger1993-py} patches to represent the middle surface of shells.
The optimized shapes are restricted by the parameterization.
In contrast,  non-parametric approaches offer a free-form design,
where  the middle surface is represented as a discrete mesh~\cite{Shimoda2014-ao}, a subdivision surface~\cite{BANDARA201862}, or the level set of a function~\cite{Kobayashi2024-gf}.
Their derivation of the shape sensitivity is typically more involved than the parametric approaches.
For simulation, the Kirchhoff-Love shell model is commonly used for analysis on the higher-order representations, such as NURBS~\cite{HUGHES20054135} or subdivision surfaces~\cite{BANDARA201862}.
When the geometry is approximated by a discrete mesh, facet-based shell elements are widely adopted in engineering for FEM simulations~\cite{Reddy2006-ur}.
Each element is considered as a small plate whose deformation follows either the Kirchhoff~\cite{Kirchhoff185051,love1888xvi} or Mindlin hypothesis~\cite{Mindlin1951-vs}.
In the graphics community, various discrete shell models~\cite{Grinspun2003,Narain2013} and solid shell elements~\cite{Montes2023} have been developed to facilitate the physical simulation of cloths, thin sheets, and geometry optimization tasks.
The primary objective of traditional form-finding of shell is to minimize bending effects under prescribed loads~\cite{Bletzinger1993-py, BLETZINGER20053438,BLETZINGER2010324,Marmo2021}, or to construct self-supported shells in a purely membrane stress state~\cite{4d-minsrf,miki2024,miki2022,Xing2024-et}.
Other methods aim to minimize compliance~\cite{JIANG2023613} or stress~\cite{ZHANG2020113036}.

	Most of the aforementioned works focus on the design of macroscopic shell structures.
	To obtain a comprehensive overview of the mechanical properties of shell lattice metamaterials, \cite{BONATTI2019301} carried out extensive finite element simulations on various TPMS-based cubic shell lattices.
	\cite{CHEN2020105288} designed several stretching-dominated shell lattices that can absorb very large energies while retaining a low density.
	\cite{MA2022110426,opencell2023} leveraged the cubic symmetry of several TPMS configurations, employing NURBS to represent the symmetric surface units and optimizing both shape and thickness to enhance isotropic stiffness.
We also note the method that integrates artificial neural networks and genetic algorithm~\cite{WANG2022115571} to predict and design the loading curves of shell lattices under compression.

\section{Preliminary}
\subsection{Shell lattice}
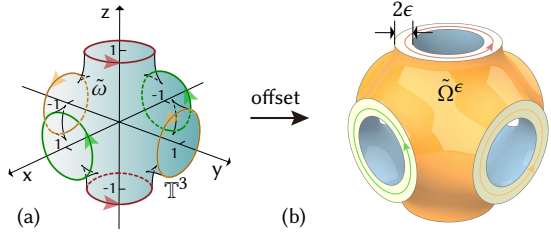
\begin{figure}[t]
	\centering
	\begin{overpic}[width=0.45\textwidth,keepaspectratio]{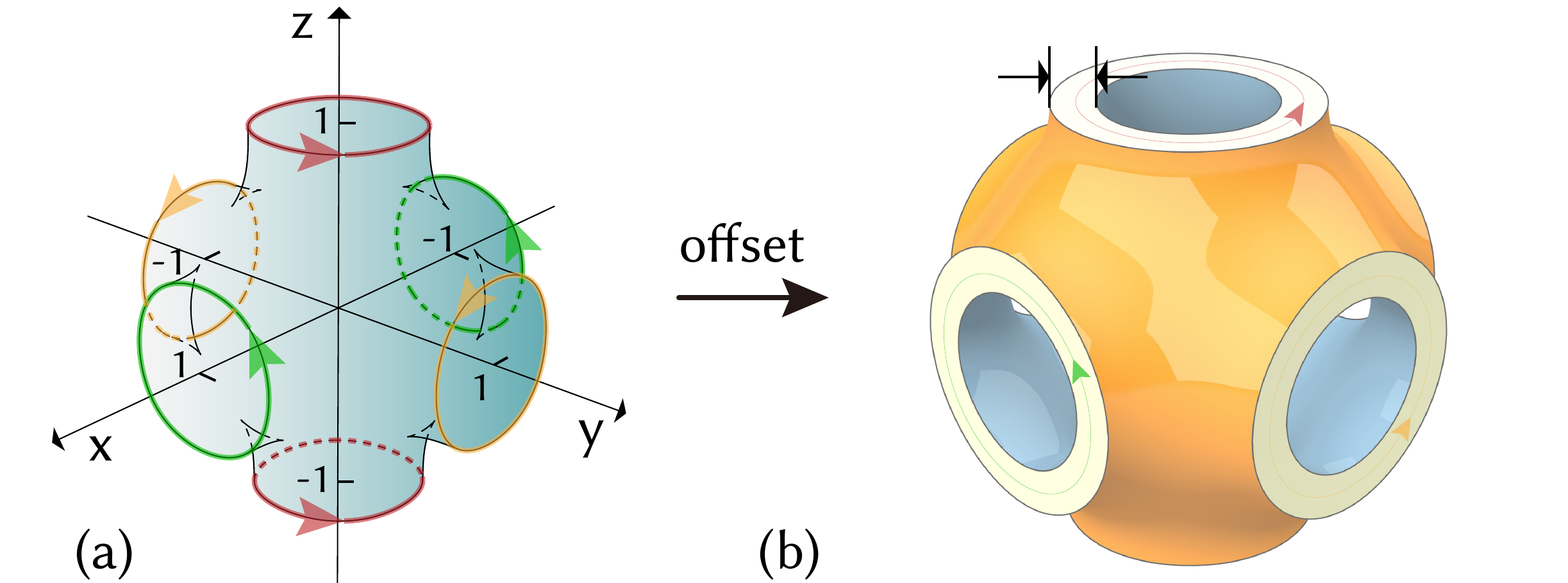}
		\put(29,5) {$\mathbb T^3$}
		\put(17,22) {$\tw$}
		\put(74,21) {$\toe$}
		\put(66.5,35) {$2\eps$}
	\end{overpic}
%	\vspace{-4mm}
	\caption{
		(a) Due to the periodic boundaries (highlighted by colored arrows), a surface $\tw$ is considered as a closed surface embedded in $\T^3$.
		(b) A shell lattice $\toe$ constructed by offsetting  $\tw$ normally on both sides. 
	}
	%    \vspace{-3mm}
	\label{fig:thick-shell}
\end{figure}

\paragraph{Construction}
The shell lattice metamaterial is constructed by thickening a smooth periodic middle surface within the RVE.
For simplicity, we consider the RVE as $Y:=[-1,1]^3$.
Due to the periodic boundaries, we identify $Y$ with a 3-torus\footnote{Note that our derivation also applies to irregular RVE, i.e., we can define the torus as $\R^3/(y_1\mathbb Z\times y_2\mathbb Z\times y_3\mathbb Z)$ for some $y_1,y_2,y_3>0$, with minor changes in proof. } $\T^3:=\R^3/(2\mathbb Z)^3$.
The middle surface is a closed surface embedded in $\T^3$, denoted as $\tw$ (Figure~\ref{fig:thick-shell}(a)).
The shell lattice, denoted as $\toe$, is constructed by offsetting $\tw$ on both sides within $\T^3$ by a distance $\eps$ (Figure~\ref{fig:thick-shell}(b)).

\paragraph{Notations}
We first introduce several necessary notations.
We denote the contravariant (covariant) components of metric tensor on $\tw$ as $a^{\a\b}$ ($a_{\a\b}$).
The second fundamental form is represented by a second-order tensor $\bs b$, with 
contravariant (covariant) components $b^{\a\b}$ ($b_{\a\b}$).
The Christoffel symbol is denoted as $\Gamma_{\a\b}^\c$.
For a vector field $\u:\tw\to\R^3$, we denote its tangential  component as $\u_{\tw}$ and the normal component as $u_3$, i.e., 
\begin{align}
	\label{eq:vec-t-n-def}
	\u_{\tw} = P\u,\quad  & u_3 = \n\cdot\u,
\end{align}
where $\n$ is the normal vector and $P=I-\n\n^\top$ denotes the projection to the tangent space.
The contravariant (covariant) components of the tangent vector $\u_{\tw}$ are denoted as $u^\a$ ($u_\a$).

\subsection{Effective stiffness}
We define the effective stiffness of lattice metamaterial under a macro-strain $\veps$  as the following energy
\begin{equation}
	\label{eq:ka-eps-matrix-form}
		\begin{aligned}
			E_\eps(\tw;{\veps})&={\veps}:{\mathbf C}_\eps:{\veps},
		\end{aligned}
\end{equation}
The symbol ${\mathbf C}_\eps$ denotes the effective elastic tensor of the shell lattice with thickness $2\eps$.
According to homogenization theory~\cite{Allaire2002}, its components are given by
\begin{equation}
	\label{eq:k-eps-wise-expr}
		\begin{aligned}
			C_\eps^{ijkl}&=\frac{1}{|Y|}\int_{\toe}\pr{\veps(\u^{ij})+\e^{ij}}:\mathbf{C}:\pr{\veps(\u^{kl})+\e^{kl}},
		\end{aligned}
\end{equation}
where $\e^{ij}=\e^i\otimes\e^j$ denotes the unit macro-strain;   $\u^{ij}$ is the solution of the following \emph{cell problem}\footnote{In the literature~\cite{elastic-texture,worst-case-stress}, the cell problem is also formulated as an equivalent differential equation.
	See the elaboration in  Section 1 of the supplementary material.
}
when ${\veps}=\e^{ij}$:
\begin{equation}
	\label{eq:min-energy}
		\begin{aligned}
			\argmin_{u\in {\bs V}_\#(\toe)} \frac{1}{|Y|}\int_{\toe}(\veps(\u) +{\varepsilon}):\mathbf C:(\veps(\u) +{\veps}).
		\end{aligned}
\end{equation}
Here, the function space is defined as $\bs V_\#(\toe):=\bs H^1(\toe)\cap\mathbf{Rig}^\perp$, where $\mathbf{Rig}^\perp$ represents the orthogonal complement  of rigid motions;
the symbols $\veps(\u)$ and ${\veps}$ denote the strain induced by $\u$ and the macro-strain, respectively;
The notation $\mathbf C$ signifies the elastic tensor of the base material constituting $\toe$, with Lam\'e coefficients $\lambda,\mu$.
Note that $E_\eps(\tw;{\veps})$ coincides with the minimal value of the objective in~\eqref{eq:min-energy}.

\section{Asymptotic analysis}
\label{sec:asym-ana}
We study the effective stiffness of shell lattices as the thickness tends to zero.
Since  it vanishes in this limit, we introduce the \emph{asymptotic directional stiffness} (ADS) to extract the leading-order behavior, defined as
\begin{align}
	\label{eq:dir-asymp-stif}
	E_A(\tw;{\veps})&:=\lim_{\eps\to 0}\frac{{E}_\eps(\tw;{\veps})}{\rho_\eps(\tw)},
\end{align}
where 
\begin{equation}
	\label{eq:rho-define}   
	\rho_\eps(\tw):=|\toe|/|Y|
\end{equation}
is the volume fraction of the solid within RVE. 

Building upon the theoretical framework of \cite{CiarletGenmem1996}, we establish the convergence theorem (Theorem~\ref{thm:asym-stif}) that provides the formulation of ADS.
The complete proof, which adapts key lemmas and techniques from \cite{CiarMem1996,CiarletGenmem1996} to handle periodic shell structures and cell problems, is presented in the supplementary material (see also \cite{hebey1996sobolev,Bernadou1994} therein).
We will discuss the fundamental properties of ADS in Section~\ref{sec:as-upper-bound}, where we establish its upper bound  and the condition to achieve this bound (Theorem~\ref{thm:asym-stif-opt-cond}).
As an application, this condition is used to justify the optimal bulk modulus of TPMS (Section~\ref{sec:opti-tpms-bulk}).

\subsection{Convergence theorem}
\label{sec:asym-stif}

We first introduce several notations before presenting our convergence theorem. 
The membrane elastic tensor $\A_\tw$ on $\tw$ is defined with components: 
\begin{equation}
	\label{eq:membrane-el-tensor-def}
	a^{\a\b\c\d}=\lambda_0a^{\a\b}a^{\c\d}+\mu\pr{a^{\a\c}a^{\b\d}+a^{\a\d}a^{\b\c}},
\end{equation}
where $\lambda_0=\frac{2\lambda\mu}{\lambda+2\mu}$ and $a^{\a\b}$ represents the metric tensor on $\tw$.
The membrane strain induced by a displacement field $\bs\eta$ on $\tw$ is expressed as the second-order tensor 
\begin{equation}
	\label{eq:gam-def}
	\gamma(\bs\eta):= \text{Sym}[\nabla \bs\eta_{\tw}]  - \eta_3 \bs b,
\end{equation}
where $\bs b$ denotes the second fundamental form and Sym denotes the symmetrization\footnote{In component form, the membrane strain is written as 
	$$ \gamma_{\a\b}(\bs\eta):=\frac{1}{2}\pr{\partial_{\a}\eta_\b+\partial_\b\eta_\a}-\Gamma_{\a\b}^\c\eta_\c-b_{\a\b}\eta_3. $$}.

\begin{remark}
	The membrane strain $\gamma(\bs\eta)$ consists of two parts.
	The first part $\text{Sym}[\nabla \bs\eta_{\tw}]$  denotes the intrinsic stretch resulted from tangent movement $\bs\eta_{\tw}$ (blue arrows in Figure~\ref{fig:ill-memstrain}).
	The second part ($-\eta_3\bs b$) denotes the stretch arising from normal offsetting (orange arrows in Figure~\ref{fig:ill-memstrain}).
	This contribution stems from the extrinsic curvature of $\tw$, which causes the normal vectors to vary spatially and become non-parallel.
	\begin{figure}[h]
		\centering
		\begin{overpic}[width=0.35\textwidth,keepaspectratio]{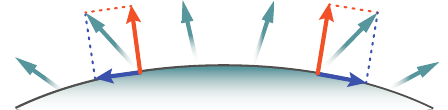}
			\put(17,24){$\bs\eta$}
			\put(20,3){$\bs\eta_{\tw}$}
			\put(30,25){$\eta_3\n$}
		\end{overpic}
%		\vspace{-2mm}
		\caption{
			Membrane strain comprises two parts: in-plane strain induced by tangent movement $\eta_{\tw}$, and the stretch resulting from normal offset $\eta_3\n$. 
		}
%		\vspace{-2mm}
		\label{fig:ill-memstrain}
	\end{figure}
	
\end{remark}

Our convergence theorem is stated as follows:

\begin{theorem}
	\label{thm:asym-stif}
	The limit in~\eqref{eq:dir-asymp-stif} is determined by
	\begin{equation}
		\label{eq:as-thm}
			E_A(\tw;\veps) = \frac{1}{|\tw|}\int_{\tw}\pr{{e}_{\tw}+\gamma(\bar\u)}:\A_{\tw}:\pr{{e}_{\tw}+\gamma(\bar \u)},
	\end{equation}
	where $e_{\tw}=P{\veps}P^\top$ is the tangential component of ${\veps}$ on $\tw$, and
	the vector field $\bar\u\in\dot{\bs{V}}_\#(\tw)$ solves the variational equation 
	\begin{equation}
		\label{eq:surf-vari-form}
			\int_{\tw}\gamma(\bar\u):\bs A_{\tw}:\gamma(\v)=-\int_{\tw}e_{\tw}: \bs A_{\tw}:\gamma(\v)\,\,\fa \v\in\dot{\bs{V}}_\#(\tw)
	\end{equation}
	in the following function space
	\begin{equation}
		\label{eq:dot-V-space-def}
		\dot{\bs{V}}_\#(\tw) := Completion\ of\ (\H^1(\tw)\,/\,\ker\gamma) \ w.r.t.\,\, norm\ |\cdot|_{\tw},
	\end{equation}
	where the norm $|\cdot|_{\tw}$ is given by
	\begin{equation}
		\label{eq:semi-norm-on-tw}
			|\v|_{\tw} := \pr{\int_{\tw}\gamma(\v):\A_{\tw}:\gamma(\v)}^{1/2}.
	\end{equation}
\end{theorem}

\begin{remark}
	The notation $\ker\gamma$ signifies the kernel of strain $\gamma$, which consists of displacements that induce zero membrane strain.
	This space is also termed as inextensional displacements in the literature.
	
\end{remark}
The detailed proof is elaborated in SM Section 4.1.
Since the solution $\bar\u$ of~\eqref{eq:surf-vari-form} is linearly dependent on the macro-strain $\veps$, the ADS admits the following compact form
\begin{equation}
	\label{eq:EA-tensor-def}
		E_A(\tw;\veps)=\veps:\mathbf C_A:\veps,
\end{equation}
where $\mathbf C_A$ is termed  as the \emph{asymptotic elastic tensor}, with components given by
\begin{equation}
	\label{eq:ca-def}
		C^{ijkl}_A :=\frac{1}{|\tw|}\int_{\tw}(\gamma(\bar\u^{ij})+\e^{ij}_{\tw}):\A_{\tw}:(\gamma(\bar\u^{kl})+\e^{kl}_{\tw}).
\end{equation}
Here, $\bar\u^{ij}$ is the solution of~\eqref{eq:surf-vari-form} for  $\veps=\e^{ij}$ and $\e^{ij}_{\tw}=P\e^{ij} P^\top$ represents the tangential  component.

\subsection{The upper bound of ADS}
\label{sec:as-upper-bound}
\begin{figure}[t]
	\centering
	\begin{overpic}[width=0.45\textwidth,keepaspectratio]{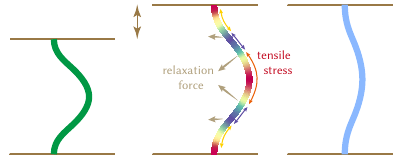}
		\put(30,33){$\veps$}
	\end{overpic}
	% \vspace{-1mm}
	\caption{
			Section view of surface deformation under a macroscopic stretch $\veps$ along the $z$-axis.
			Left: undeformed surface. 
			Middle: the RVE and surface undergo affine deformation, resulting in strain energy ($E_M$).
			The projected in-plane strain varies across the surface, leading to  non-uniform tensile stress whose divergence induces a relaxation force.
			Right: the surface relaxes under this force to reach a lower energy state ($E_A$).
	}
	% \vspace{-2mm}
	\label{fig:ill-relax}
\end{figure}

By substituting the variational equation~\eqref{eq:surf-vari-form} into~\eqref{eq:as-thm} with $\v=\bar\u$, the ADS can be expressed as (Appendix~\ref{app:compact-form-EM})
\begin{equation}
	\label{eq:EA-EM-neg}
		E_A(\tw;\veps) = E_M(\tw;\veps) - \frac{1}{|\tw|}\int_{\tw}\gamma(\bar\u):\A_{\tw}:\gamma(\bar\u).
\end{equation}
We define $E_M(\tw;\veps)$ as the \emph{homogeneous membrane energy}, given by
\begin{equation}
	\label{eq:EM-defin}
		E_M(\tw;\veps):=\veps:\overline{\mathbb{P}}:\veps,
\end{equation}
where
$
\overline{\mathbb{P}}:=\frac{1}{|{\tw}|}\int_{\tw}\mathbb{ P}
$
represents the average of the fourth-order tensor $\mathbb P$ over $\tw$, defined as
\begin{equation}
	\label{eq:P-form-tensor}
	[\mathbb P]^{ijkl}=\lambda_0 P^{ij}P^{kl}+\mu\pr{P^{il}P^{jk}+P^{ik}P^{jl}}.
\end{equation}
The detailed derivation is provided in Appendix~\ref{app:compact-form-EM}.
Since $\A_{\tw}$ is a positive definite tensor,  we have
$
	\frac{1}{|\tw|}\int_{\tw}\gamma(\bar\u):\A_{\tw}:\gamma(\bar\u)\ge 0
$
with equality holding if and only if $\bar \u=0$.
This leads to the following theorem that establishes the upper bound of $E_A(\tw;\veps)$.
\begin{theorem}
	\label{thm:EA-upper-limit}
	The asymptotic directional stiffness satisfies
	\begin{equation}
		E_A(\tw;\veps) \le E_M(\tw;\veps),
	\end{equation}
	and the equality holds if and only if~\eqref{eq:surf-vari-form} admits a zero solution.
\end{theorem}
\begin{remark}
	\label{rem:intui-upp-bd}
	The intuition behind this theorem is as follows: 
		When the surface is subjected to a macro-strain $\veps$, the energy $E_M$ corresponds to the strain energy under purely affine deformation. However, due to the presence of void space around the surface within the RVE, the surface can further relax, reducing the energy to the final value $E_A$.
		If no relaxation occurs  (i.e., $\bar\u=\bs 0$), then the final stable energy $E_A$ equals $E_M$.
		See Figure~\ref{fig:ill-relax} for an illustration.
	
\end{remark}

\begin{remark}
	\label{thm:eigen-sum-of-P}
	Although $E_M$ still depends on $\veps$ as  does $E_A$, it has a simpler structure since the `matrix' $\overline{\mathbb P}$ of the quadratic form~\eqref{eq:EM-defin} only depends on the normal vector distribution of $\tw$, facilitating analysis.
	For example, we know the eigenvalues of tensor $\overline{\mathbb P}$ sum to $2\lambda_0+6\mu$, namely, 
	\begin{equation}
		\frac{1}{6}\sum_{i=1}^6 \lambda_i(E_M)=\frac{1}{3}\lambda_0+\mu,
	\end{equation}
	where $\lambda_i$ denotes the $i$-th eigenvalue.
	See the proof in Appendix~\ref{app:eigen-sum-of-P}.
	This result implies that while $E_M$ varies with $\veps$, its average eigenvalue is constant.
\end{remark}

\subsection{Achieving the upper bound of ADS}
\label{sec:condi-upper-bound}
One common objective in lattice metamaterial design is stiffness optimization~\cite{opencell2023,Liu2024-multilayer}, motivating our search for surfaces achieving the upper bound $E_M$.
While Theorem~\ref{thm:EA-upper-limit} characterizes the optimal condition, its direct application as a design criterion is impractical due to the need to solve variational equation~\eqref{eq:surf-vari-form}.
Therefore, we introduce an equivalent necessary and sufficient condition:  
\begin{theorem}
	\label{thm:asym-stif-opt-cond}
	The asymptotic directional stiffness satisfies
	\begin{equation}
		E_A(\tw;\veps) \le E_M(\tw;\veps),
	\end{equation}
	and the equality holds if and only if the following relations hold over $\tw$:
	\begin{equation}
		\label{eq:opt-asym-stiff-condi}
		\begin{aligned}
			\pr{2\pr{\lambda_0+\mu}\bs b+ 4\mu H P} \veps\cdot \n=\mathbf 0\\
			2\pr{\lambda_0 H P+\mu \bs b}:\veps  =0,
		\end{aligned}  
	\end{equation}
	where the second fundamental form is considered as a tensor in $\T^3$ during the contraction\footnote{This is feasible as $\tw$ is smoothly embedded in $\T^3$.}.
\end{theorem}
\begin{remark}
	\label{rem:relax-force}
		The tangent vector and scalar on the LHS of~\eqref{eq:opt-asym-stiff-condi} are exactly the tangential and normal components of the relaxation force acting on the surface, as shown in Figure~\ref{fig:ill-relax}.
\end{remark}

The complete proof of Theorem~\ref{thm:asym-stif-opt-cond} is provided in SM Section 4.2. 
Note that this condition only involves the geometry of the middle surface without solving any equation. 
It suggests that constructing a surface satisfying \eqref{eq:opt-asym-stiff-condi} universally for all $\veps$ would yield a optimal (in the sense of stiffness) middle surface for designing shell lattice. 
Unfortunately, we have the following theorem
\begin{theorem}
	\label{thm:no-opt-surfaces}
	The only smooth surface such that \eqref{eq:opt-asym-stiff-condi} holds everywhere for every $\veps$ is plane.
\end{theorem}
\begin{proof}
	Suppose $\tw$ satisfies \eqref{eq:opt-asym-stiff-condi} at every point for every $\veps$, then we have 
	\begin{equation}
		\label{eq:hp-b-eq-zero}
		2\pr{\lambda_0 HP+\mu\bs b} =\mathbf 0.
	\end{equation}
	It follows that
	\begin{equation}
		\tr 2\pr{\lambda_0 HP+\mu\bs b}=\pr{4\lambda_0+2\mu}H=0,
	\end{equation}
	where we have utilized~\eqref{eq:trace-of-P} and the identity $\tr\bs b=2H$.
	The above equation implies zero mean curvature: $H=0$.
	Hence we have $\bs b =\mathbf 0$ according to~\eqref{eq:hp-b-eq-zero}.
	This means $\tw$ is a plane in $\T^3$.
\end{proof}

Nevertheless, non-trivial surfaces that achieve the upper bound exist if we only consider certain strain, as we will show in the next section.

\section{The optimal asymptotic bulk modulus of TPMS}
\label{sec:opti-tpms-bulk}
We define the \emph{asymptotic bulk modulus} (ABM) as the ADS under the hydrostatic strain ($\veps=\frac{1}{3}I$), i.e., 
\begin{equation}
	K_A(\tw) := E_A(\tw;\frac{1}{3}I).
\end{equation}
Theorem~\ref{thm:EA-upper-limit} gives the following upper bound on ABM:
\begin{equation}
	K_A(\tw) \le E_M(\tw;\frac{1}{3}I).
\end{equation}
Surprisingly, we have the following proposition
\begin{proposition}
	\label{thm:EM-hydrostatic-strain}
	The homogeneous membrane energy under hydrostatic loading is constant:
	\begin{equation}
		E_M(\tw;\frac{1}{3}I) = \frac{4}{9}\pr{\lambda_0+\mu}.
	\end{equation}
\end{proposition}
\begin{proof}
	Direct computation shows  
	\begin{equation}
		\begin{aligned}
				E_M(\tw;\frac{1}{3}I)&=\frac{1}{9|\tw|}\int_{\tw}I:{\mathbb P}:I=\frac{1}{9|\tw|}\int_{\tw}\lambda_0 \pr{\tr P}^2 +2\mu\tr PP^\top\\&=\frac{1}{9|\tw|}\int_{\tw}4\pr{\lambda_0+\mu}=\frac{4}{9}\pr{\lambda_0+\mu},
		\end{aligned}
	\end{equation}
	where we have used~\eqref{eq:P-form-tensor} and~\eqref{eq:trace-of-P}.
\end{proof}

The interesting problem is which kind of surface attains this upper bound.
If we substitute $\veps=\frac{1}{3}I$ into LHS of~\eqref{eq:opt-asym-stiff-condi}, we  get
\begin{equation}
	\label{eq:relax-force-bulk}
	\begin{aligned}
		\pr{2\pr{\lambda_0+\mu}\bs b+ 4\mu H P} \veps\cdot \n&\to \bs 0\\
		2\pr{\lambda_0 HP+\mu\bs b}:\veps&\to\frac{4}{3}\pr{\lambda_0 +\mu}H.
	\end{aligned}
\end{equation}
Here, the first relation holds because $\bs b$ and $P$ are both orthogonal to the normal vector; the second  utilized~\eqref{eq:trace-of-P} and the identity $\tr\bs b=2H$.
Therefore, the optimality condition~\eqref{eq:opt-asym-stiff-condi} reduces to $H=0$ in this case.
Consequently, we have derived the following theorem.
\begin{theorem}
	\label{thm:tpms-opt-bulk}
	The asymptotic bulk modulus satisfies
	\begin{equation}
		\label{eq:KA-upper-bound}
		K_A(\tw)\le \frac{4}{9}\pr{\lambda_0+\mu},
	\end{equation}
	and the equality holds if and only if $\tw$ is TPMS.
\end{theorem}

\begin{remark}
		As discussed in Remark~\ref{rem:relax-force}, equation~\eqref{eq:relax-force-bulk} implies that the relaxation force under hydrostatic loading is a normal force proportional to the mean curvature. 
		Therefore, if the mean curvature vanishes, no relaxation occurs, and the ADS reaches its upper bound.
\end{remark}

\paragraph{Relation to Hashin-Shtrikman bound}
This theorem explains the  exceptionally high bulk modulus observed in TPMS shell lattice metamaterials~\citep{pnas-shell-limit,CHEN2019108109,D0MH01408G}, which approaches the Hashin-Shtrikman (HS) upper bound~\citep{HASHIN1963127}.
Here, we show how the theorem relates to this bound.
According to the composite material theory, the HS upper bound for bulk modulus is given by  
\begin{equation}
	K^{HS}(\rho):=\frac{4\rho\mu(3\lambda+2\mu)}{12\mu+(1-\rho)(9\lambda+6\mu)},
\end{equation}
where $\rho$ is the volume fraction of the material within the unit cell. 
When this bound applies to our problem, we immediately get  
\begin{equation}
	\label{eq:bulk-hs-bound}
	E_\eps(\tw;\frac{1}{3}\mathbf I)= K(\toe)\le K^{HS}(\rho(\toe)).
\end{equation}
If we divide $\rho(\toe)$ on both sides of~\eqref{eq:bulk-hs-bound} and let $\epsilon\to 0$, we have
\begin{equation}
	\label{eq:abm-hs-bound}
	\limsup_{\epsilon\to 0}\frac{E_\eps(\tw;\frac{1}{3}\mathbf I)}{\rho(\toe)}\le\lim_{\epsilon\to 0}\frac{K^{HS}(\rho(\toe))}{\rho(\toe)}=\frac{\itd K^{HS}}{\itd\rho}(0)=\frac{4}{9}(\lambda_0+\mu),
\end{equation}
which is exactly the upper bound in \eqref{eq:KA-upper-bound}.
However, this derivation necessitates that the lattice metamaterial exhibits cubic symmetry; otherwise, the HS bound does not generally hold. Our approach relaxes this constraint, provided that the shell lattice has a smooth, periodic middle surface.
As a result, this theorem implies that any cubic TPMS shell lattice  consistently exhibits a bulk modulus close to HS upper bound (at low volume fractions).
Even when the surface lacks cubic symmetry, this theorem remains valid, though the HS bound may no longer apply.
We stress that~\eqref{eq:KA-upper-bound} holds  even when the surface is not closed or the thickness is not uniform (SM Section 4.3).

\section{Discretization and optimization}
\label{sec:disc}
The optimality condition~\eqref{eq:opt-asym-stiff-condi} in Theorem~\ref{thm:asym-stif-opt-cond} offers a guiding principle for designing surfaces with extreme stiffness.
However, constructing a surface that satisfies~\eqref{eq:opt-asym-stiff-condi} everywhere for an arbitrary strain is not straightforward, except in the special case of bulk modulus.  
	Moreover, the upper bound itself depends on the surface geometry, potentially limiting the achievable stiffness.
To address these issues, we introduce a numerical framework for   optimizing ADS.
We first illustrate our discretization scheme for evaluating ADS (Section~\ref{sec:ads-discre}) before presenting the optimization algorithm (Section~\ref{sec:shape-opt}).

\subsection{Discretization of ADS}
\label{sec:ads-discre}
We approximate the middle surface $\tw$ by a triangular mesh, denoted as $\mathcal M$.
The sets of faces, edges, and vertices of $\mathcal M$ are denoted as $\mathcal F$, $\mathcal E$, and $\mathcal V$, respectively.

\begin{figure}[t]
	\centering
	\begin{overpic}[width=0.45\textwidth,keepaspectratio]{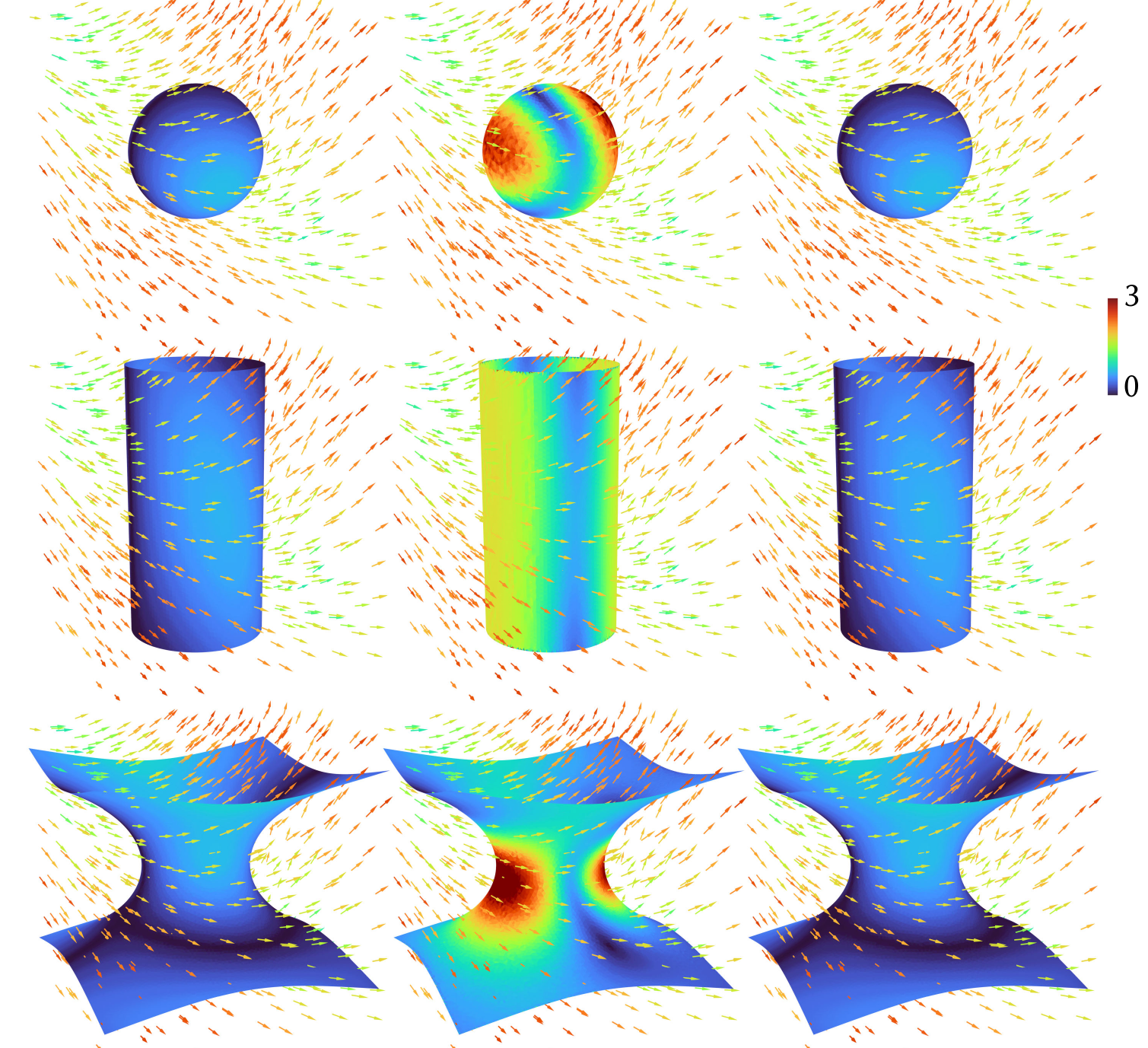}
	\end{overpic}
%	\vspace{-5mm}
	\caption{
		Approximation error of the membrane strain on elliptic, parabolic, and hyperbolic surfaces under a smooth vector field $(\cos(x),\sin(xy),\sin(x+z))$  (visualized with colored arrows).
		Left : the plane stress element~\eqref{eq:plane-element}.
		Middle : scheme~\eqref{eq:second-gam-approx}.
		Right : scheme~\eqref{eq:gam-disc-final}.
	}
%	\vspace{-1mm}
	\label{fig:accu-gamma}
\end{figure}

\subsubsection{The membrane strain}
\label{sec:the-mem-strain}
\begin{figure}[t]
	\centering
	\begin{overpic}[width=0.45\textwidth,keepaspectratio]{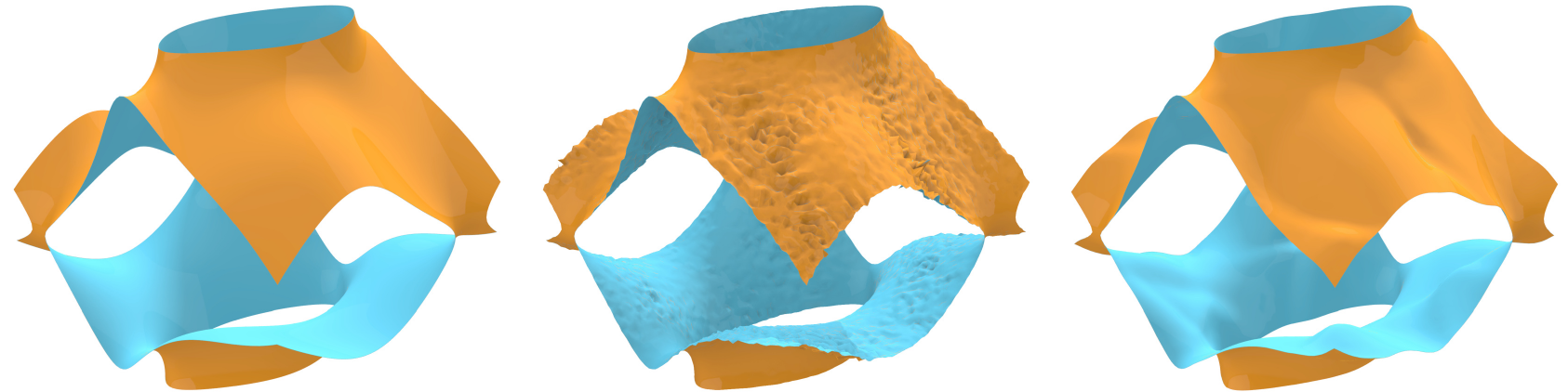}
	\end{overpic}
%	\vspace{-5mm}
	\caption{
		Deformation solved from~\eqref{eq:surf-vari-form} under hydrostatic loading ($\veps=\frac{1}{3}I$) with different elements.
		Left : the input surface.
		Middle :  plane stress element~\eqref{eq:plane-element}.
		Right : scheme~\eqref{eq:gam-disc-final}.
	}
	%    \vspace{-4mm}
	\label{fig:plate-singular}
\end{figure}

By definition, $\gamma(\u)$ is the change of metric when $\tw$ deforms along a vector field $\u$.
Thus, a direct approximation computes the strain of each triangular face
	$\bigtriangleup \bs{x}_1\bs{x}_2\bs{x}_3$ under vertex displacements $\u_i$ ($i=1,2,3$) as
\begin{equation}
	\label{eq:plane-element}
	\gamma(\u)\approx \text{Sym}[\nabla(P_f\u)],
\end{equation}
where $P_f$ denotes the projection onto face $f$.
Essentially, this is a plane stress element with no bending effects taken into account.
The approximation, however, results in a highly fluctuating displacement solution of~\eqref{eq:surf-vari-form}, which is physically implausible (Figure~\ref{fig:plate-singular} middle).
As suggested by the definition~\eqref{eq:gam-def} of $\gamma(\u)$, 
this issue may stem from the lack of curvature information in the surface.
To address this, we attempt to approximate the second fundamental form to better capture the surface’s geometry. 

\paragraph{Approximating the second fundamental form}
We apply finite difference scheme~\cite{Rusinkiewicz2004} to approximate the second fundamental form.
According to the definition, we have 
\begin{equation}
	\label{eq:b-eval-edge}
	\bs{b}(\bs{l}_{ij}, \bs{l}_{ij})=-(\n_j-\n_i)\cdot \bs{l}_{ij},
\end{equation}
where $\bs n_i$ denotes the normal at $i$-th vertex, computed using the angle based weighting scheme, and $\bs{l}_{ij}=\bs v_j-\bs v_i$ denotes the edge vector.
On each face $f\in\mathcal F$, we approximate $\bs b$ by a symmetric matrix $\bs b_f\in\mathbb R^{2\times2}$ with three independent components determined from its evaluation~\eqref{eq:b-eval-edge} over the edge vectors of the face.
On the other hand, the normal displacement $u_3$ within the face is approximated by linear interpolation from the displacements at the vertices:
%\vspace{-1mm}
\begin{equation}
	u_3 =\sum_{i=1}^3 \phi_i \u_i\cdot\n_i,
\end{equation}
%\vspace{-1mm}
where $\phi_i$ denotes the barycentric coordinate.
Thus, we refine our approximation of the membrane strain as follows:
%\vspace{-1mm}
\begin{equation}
	\label{eq:second-gam-approx}
	\gamma(\bs u) \big|_f \approx \text{Sym} [\nabla (P_f\bs u)] - u_3\bs b_f.
\end{equation}
%\vspace{-1mm}
Unfortunately, this formulation exhibits poor accuracy in experiments (Figure~\ref{fig:accu-gamma} middle).
We find that this significant error actually stems from the symmetric gradient in~\eqref{eq:second-gam-approx}, instead of the second fundamental form.

\begin{wrapfigure}{r}{0.5\columnwidth}
	\hspace{-20pt}
	\vspace{-15pt}
	% \centering
	\begin{overpic}[width=0.99\linewidth]{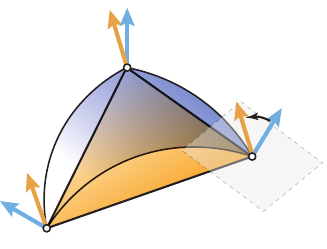}
		\put(23,68){$\n_f$}
		\put(43,70){$\n_j$}
		\put(37,47){$\v_j$}
		\put(10,0){$\v_k$}
		\put(75,20){$\v_i$}
		\put(88,36){$\n_i$}
		\put(77,44){$\bs R_i^f$}
	\end{overpic}
	\caption{
	}
	\label{fig:vert-proj}
\end{wrapfigure}

\paragraph{Correction of tangent strain}
To illustrate this issue, we consider the approximation of a curved triangular surface patch  by a planar triangular face (Figure~\ref{fig:vert-proj}).
We assume that the vertex positions and normals coincide (blue arrows in Figure~\ref{fig:vert-proj}), and that the curvature is small relative to the dimensions of both patches.
In this case, the  tangent space of the surface can be well approximated by the plane of the  triangle, along with the tangent vectors defined on it.
However, this approximation is in the sense of $L^2$ norm.
Therefore, it is not accurate to approximate the gradient of the tangent vectors by differentiating this $L^2$ approximation, as~\eqref{eq:second-gam-approx} does.
To address this issue, we must account for the curvature of the surface, which is characterized by the non-parallel vertex normals.
This observation inspires us to use per-vertex projection instead of a uniform projection ($P_f$) to extract the tangential  component. 
In practice, we first rotate the vertex tangent space to align with the plane of the face before projecting the displacement (Figure~\ref{fig:vert-proj}).
Then we approximate the tangent vector within the face as the linear interpolation of the projected ones at the vertices:
%\vspace{-1mm}
\begin{equation}
	\label{eq:vertex-proj}
	\u_{\tau} = \sum_{i=1}^3\phi_i P_f R_{i}^{\tiny f} \u_i,
\end{equation}
%\vspace{-1mm}
where $\phi_i$ denotes barycentric coordinate, and $R^f_i$ is a rotation matrix aligning vertex normal $\n_i$ to face normal $\n_f$ (yellow arrows in Figure~\ref{fig:vert-proj}).
Finally, our improved approximation is given by 
\begin{equation}
	\label{eq:gam-disc-final}
	\gamma(\u)|_f \approx \text{Sym}[\nabla \u_\tau] -u_3\bs{b}_f.
\end{equation}
The numerical experiments demonstrate a significant improvement in accuracy after this correction (Figure~\ref{fig:accu-gamma} right).

Although the plane stress element~\eqref{eq:plane-element} and our scheme exhibit comparable accuracy in discretizing the strain w.r.t. a smooth vector field (Figure~\ref{fig:accu-gamma} left and right), our approach yields a more convincing solution to~\eqref{eq:surf-vari-form} (Figure~\ref{fig:plate-singular} right).
The accuracy and convergence rate of this scheme are further investigated  in Section~\ref{sec:num-problem}.

\subsubsection{Solve for ADS}
\label{sec:the-vari-eq}

With the discretized membrane strain, the variational equation~\eqref{eq:surf-vari-form} is converted into the following linear system (Appendix~\ref{eq:assemb-stif-matvec}):
\begin{equation}
	\label{eq:surf-equation-stif}
	\mathbf{K}^h\u^h = \mathbf{f}^h,
\end{equation}
where $\mathbf{K}^h$ denotes the global stiffness matrix, and $\mathbf{f}^h$ is the force vector assembled from the macro-strain $\veps$.
The vector $\u^h$ consists of the vertex displacements that we solve for.
We denote the force vector as $\mathbf{f}^h_{ij}$ when $\veps =\e^{ij}$, with the corresponding solution of~\eqref{eq:surf-equation-stif} denoted as $\u^h_{ij}$.
Finally, the components of asymptotic elastic tensor are computed as:
%\vspace{-2mm}
\begin{equation}
	\label{eq:disc-ca-form}
	C_A^{ijkl} = E_M^h(\e^{ij},\e^{kl}) - \frac{\u^h_{ij}\cdot\mathbf f_{kl}^h}{Area(\mathcal M)},
\end{equation}
%\vspace{-1mm}
where
%\vspace{-1mm}
\begin{equation}
	E_M^h(\e^{ij},\e^{kl}):=\frac{1}{Area(\mathcal M)}\sum_{f\in\mathcal F}\int_{f} \A_{\tw}:P_f\e^{ij}P_f^\top:P_f\e^{kl}P_f^\top.
\end{equation}
%\vspace{-2mm}
We have used the variational equation~\eqref{eq:surf-vari-form} to derive~\eqref{eq:disc-ca-form}, similar to~\eqref{eq:EA-EM-neg}.

\subsection{Shape optimization}
\label{sec:shape-opt}
As shown in the definition of membrane strain, computing the second fundamental form and the Christoffel symbols requires second-order derivatives of the surface's embedding.
Although it can be approximated piecewise (Section~\ref{sec:the-mem-strain}), the shape sensitivity involves third-order derivatives, which makes computation on a triangular mesh challenging.
While higher-order geometric representations such as NURBS or Bézier patches can address this issue, their generation and optimization on a periodic surface with complex topology is non-trivial.
To overcome this challenge, we leverage the variational equation~\eqref{eq:surf-vari-form} to eliminate the higher-order terms, resulting in a feasible formulation that is easier to discretize (Section~\ref{sec:sens-ads}).
This sensitivity is then used to derive a gradient-based flow for optimizing ADS (Section~\ref{sec:geo-opt}).

\subsubsection{Sensitivity of ADS}
\label{sec:sens-ads}
As tangential  movement does not change the surface shape, we assume that $\tw$ evolves under a normal velocity $v_n:\tw\to\R$.
The time derivative of ADS, which by formulation~\eqref{eq:as-thm} is a ratio of an integral (denoted as $I_a^s$) to the surface area (denoted as $A$), takes the following form
\begin{equation}
	\label{eq:common-ad-rate}
	\dt E_A=\frac{\dt{I}^s_a}{A}-\frac{I_a^s \dt{A}}{A^2}.
\end{equation}
The change of area is straightforward to compute, given by
\begin{equation}
	\label{eq:area-rate}
	\dt{A}:=-\int_{\tw}2v_n H.
\end{equation}
The discretization of this integral is discussed in Appendix~\ref{app:sens}.
Then it remains to compute the time derivative of the integral $\dt{I}_a^s$. 

The vanilla derivation gives the following time derivative (SM Section 8):

%\vspace{-1mm}
\begin{proposition}
	\label{prop:ads-sens}
	The time derivative $\dt{I}_a^s$ is given by
%	\vspace{-1mm}
	\begin{equation}
		\label{eq:sens-ads-form-raw}
		\begin{aligned}
				\dt{I}_a^s&=\int_{\tw}2 v_n\pr{\gamma(\bar\u)+e_{\tw}}:(\bs{B}-H\bs{A}_{\tw}):\pr{\gamma(\bar\u)+e_{\tw}}\\
				&+\int_{\tw}2\pr{\gamma(\bar\u)+e_{\tw}}:\A_{\tw}:\xi(\bar\u),
		\end{aligned}
	\end{equation}
%	\vspace{-1mm}
	where $\bs B$ is a fourth-order tensor defined as
	\begin{equation}
		\label{eq:bijkl-def}
		B^{ijkl} := \lambda_0(b^{\a\b}a^{\c\d}+a^{\a\b}b^{\c\d})+\mu(b^{\a\c}a^{\b\d}+a^{\a\c}b^{\b\d}+b^{\a\d}a^{\b\c}+a^{\a\d}b^{\b\c}),
	\end{equation}
	and $\xi(\bar\u)$ is a second-order tensor, defined as
	\begin{equation}
		\label{eq:ads-sens-xi-def}
		\begin{aligned}
			\xi(\u)&:=v_n\bar\u_{\tw}\cdot\nabla\bs{b}+2\text{Sym}[\bs b\bar\u_{\tw} \nabla v_n^\top]-(\nabla v_n\cdot\bar\u_{\tw})\bs b\\
			&-\bar{u}_3\pr{\nabla^2 v_n-v_n\bs c}+2\text{Sym}[P\veps\n\nabla v_n^\top].
		\end{aligned}
	\end{equation}
	where $\bs c=\bs b\cdot\bs b$ denotes the third fundamental form.
\end{proposition}

%\vspace{-1mm}
A challenge with this formulation is the presence of the covariant derivative of the second fundamental form ($\nabla\bs b$) and the second-order derivative ($\nabla^2 v_n$), both of which are difficult to evaluate on a discrete mesh. 
While the existing work~\cite{Rusinkiewicz2004} proposed an approximation for the curvature derivatives on discrete meshes, we introduce a more refined approach that bypasses this approximation.

\begin{figure}[t]
	\centering
	\begin{overpic}[width=0.48\textwidth,keepaspectratio]{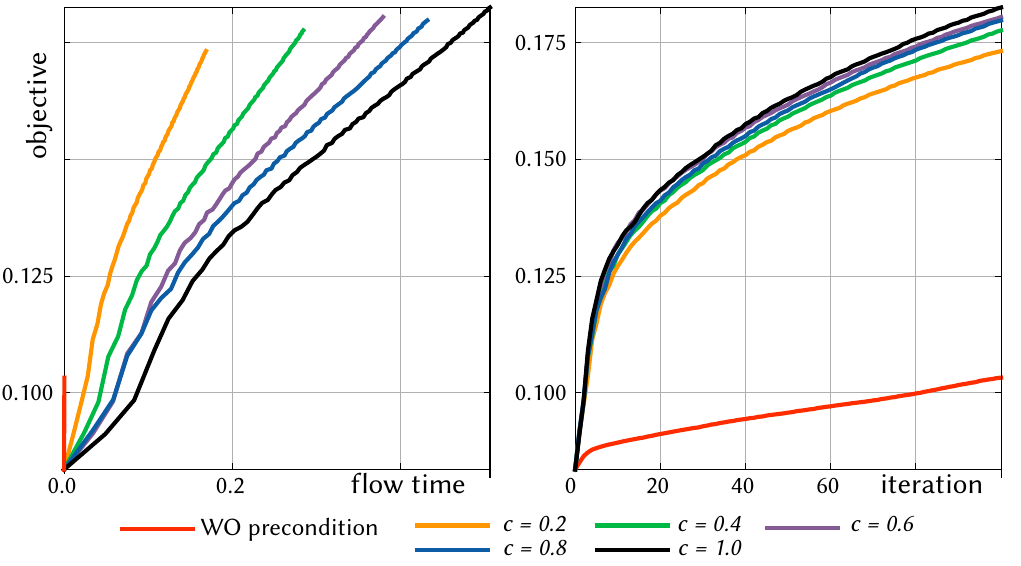}
	\end{overpic}
	%	\vspace{-7mm}
	\caption{
		Maximizing ADS under strain $\e^{33}$ without precondition and with different precondition strength parameter $c$.
		The left and right plots show the objective change w.r.t. the flow time (parameter $t$ in~\eqref{eq:flow-eq})  and the iteration, respectively.
		Although the objective evolves fastest w.r.t. flow time without preconditioning (left), the total clock time required to reach the same objective value depends on the number of iterations, which is significantly reduced  after preconditioning (right).
	}
	%    \vspace{-2mm}
	\label{fig:grad-precondition}
\end{figure}

\paragraph{Eliminate the higher-order derivatives}
To handle this problem, we notice that both terms can be written as
\begin{align}
	v_n\bar{\u}_{\tw}\cdot\nabla\bs{b}&=\nabla(v_n\bs{b}\bar\u_{\tw})-\bs{b}\nabla(\v_n\bar\u_{\tw})\\
	\bar u_3\nabla^2 v_n &= \nabla (\bar u_3\nabla v_n)-\nabla\bar u_3\nabla v_n^\top.
\end{align}
The benefit is that we can leverage the variational equation~\eqref{eq:surf-vari-form} to eliminate both terms by substituting the test function with $\v=v_n\bs{b}\bar\u_{\tw}$ (resp. $\bar u_3\nabla v_n$).
Similarly, the term $(\nabla v_n\cdot\bar\u_{\tw})\bs b$ is  eliminated by letting $v_3=(\nabla v_n\cdot\bar\u_{\tw})$ and $\v_{\tw}=\bs 0$.
\begin{remark}
		A potential issue is that the substituted test functions may not belong to the space $\dot{\bs V}_\#(\tw)$, particularly those involving $\bar\u_{\tw}$. 
		As a result, the above approach implicitly assumes sufficient regularity of both the solution $\bar\u$ and the normal velocity $v_n$.
\end{remark}
After collecting the remaining terms, we can replace $\xi(\bar\u)$ in~\eqref{eq:sens-ads-form-raw} with 
\begin{equation}
	\label{eq:impro-sens-zeta}
	\zeta(\bar\u):=\bs{b}\bar\u_{\tw}\nabla v_n^\top -v_n\bs b\nabla\bar\u_{\tw}+\nabla \bar u_3\nabla v_n^\top +\bar u_3 v_n\bs c+2P\veps\bs n\nabla v_n^\top,
\end{equation}
where we omit the symmetrization operation for brevity, as its  contracted tensor in~\eqref{eq:sens-ads-form-raw} is symmetric.
Thus, we propose the following improved sensitivity formulation
\begin{proposition}
	\label{prop:ads-sens-impro}
	The time derivative $\dt{I}_a^s$ is evaluated as
	\begin{equation}
		\label{eq:sens-ads-form-raw-impro}
		\begin{aligned}
				\dt{I}_a^s&=\int_{\tw}2 v_n\pr{\gamma(\bar\u)+e_{\tw}}:(\bs{B}-H\bs{A}_{\tw}):\pr{\gamma(\bar\u)+e_{\tw}}\\
				&+\int_{\tw}2\pr{\gamma(\bar\u)+e_{\tw}}:\A_{\tw}:\zeta(\bar\u),
		\end{aligned}
	\end{equation}
	where $\bs B$ is defined in~\eqref{eq:bijkl-def} and $\zeta(\bar\u)$ is defined in~\eqref{eq:impro-sens-zeta}.
\end{proposition}
This new formulation contains no higher-order derivatives, resulting in a suitable discretization on a triangle mesh (Appendix~\ref{sec:ads-discre}).

\subsubsection{Sensitivity of asymptotic elastic tensor}
Since ADS is a quadratic form, we can derive the sensitivity of its associated matrix, the asymptotic elastic tensor $\mathbf C_A$, based on the sensitivity of its evaluation on strains.
Alternatively, one may directly compute the time derivative of the integral in the definition~\eqref{eq:ca-def}, following a similar procedure as in Section~\ref{sec:sens-ads}.
Both approaches lead to the same result, summarized below:
\begin{proposition}
	The time derivative of asymptotic elastic tensor $\mathbf C_A$ is given by
%	\vspace{-1mm}
	\begin{equation}
		\label{eq:ca-ijkl-dt}
		\begin{aligned}
			\dt{C}^{ijkl}_A &= \frac{\dt{I}_{ijkl}^s}{A} -\frac{I_{ijkl}^s\dt{A}}{A^2},\\
		\end{aligned}
	\end{equation}
	where $I^s_{ijkl}$ denotes the integral in~\eqref{eq:ca-def}, respectively; and their time  derivatives are
	\begin{align}
		\dt{I}^s_{ijkl}&=\int_{\tw} 2v_n\pr{\gamma(\bar\u^{ij})+\e^{ij}_{\tw}}:\pr{\bs B-H\A_{\tw}}:\pr{\gamma(\bar\u^{kl})+\e^{kl}_{\tw}}\\
		&+\int_{\tw}\pr{\gamma(\bar\u^{ij})+\e^{ij}_{\tw}}:\A_{\tw}:\zeta(\bar\u^{kl})\\&+\int_{\tw}\pr{\gamma(\bar\u^{kl})+\e^{kl}_{\tw}}:\A_{\tw}:\zeta(\bar\u^{ij}).
	\end{align}
\end{proposition}
%\vspace{-1mm}
The derivation follows the same spirit in Section~\ref{sec:sens-ads} with minor changes in the integrand, and is therefore omitted here.

\subsubsection{Optimization}
\label{sec:geo-opt}

According to~\eqref{eq:area-rate},~\eqref{eq:common-ad-rate}, and~\eqref{eq:sens-ads-form-raw}, the time derivative $\dt{I}_A^s$ is a linear functional of the normal velocity $v_n$.
Hence, we can define the $L^2$ gradient $\grad I_A^s$ such that
%\vspace{-1mm}
\begin{equation}
	\label{eq:l2-graident}
	\dt{I}_A^s(v_n) = \inn{\grad I_A^s, v_n}_{L^2(\tw)},
\end{equation}
%\vspace{-1mm}
and then performing a gradient descent flow as
%\vspace{-1mm}
\begin{equation}
	\label{eq:flow-eq}
	\frac{d\tw}{dt} = -\grad I_A^s. %\qquad\qquad\qquad
\end{equation}
\begin{wrapfigure}{r}{0.25\columnwidth}
	\hspace{-8pt}
	\begin{overpic}[width=0.99\linewidth]{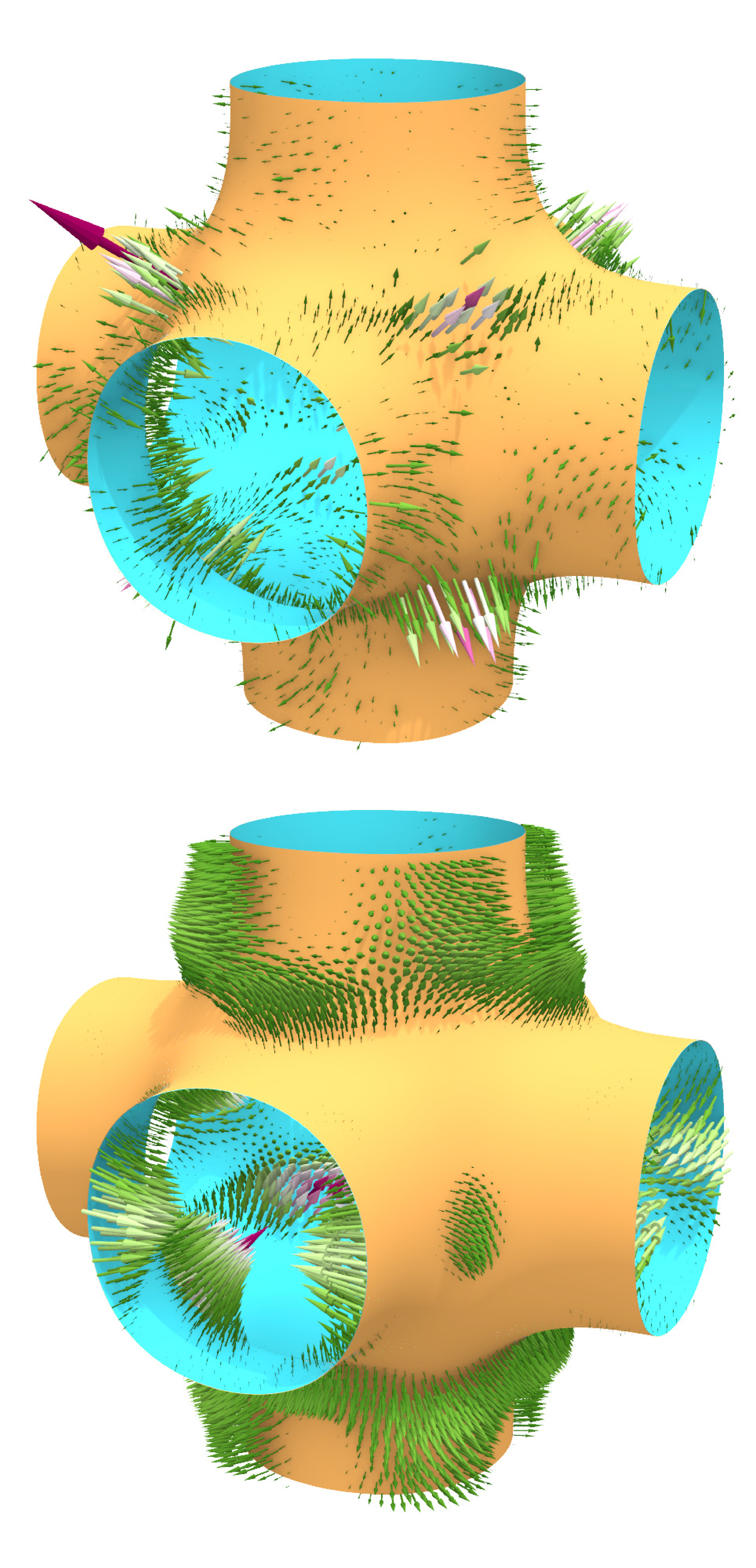}
	\end{overpic}
	\vspace{-10pt}
	\hspace{-8pt}
	\caption{
	}
	\label{fig:grad-prec-arrow}
\end{wrapfigure}
However, the Armijo condition during line search forces a significant reduction in the time step when using the $L^2$ gradient, resulting in very small time steps even if a large initial step is specified, thereby slowing down the optimization (Figure~\ref{fig:grad-precondition} right).
Similar issues are also noted in~\cite{Yu2021-crv,Yu2021-srf}, 
where they are addressed by using the gradient derived from fractional Sobolev inner product that matches the energy objective.
Here, we resort to another approach.
We precondition this gradient with the Laplacian operator by solving for a modified gradient $d_n$ from the following equation
%\vspace{-1mm}
\begin{equation}
	(Id-c \Delta ) d_n = - (c+1)\grad I_A,\qquad\qquad\qquad\qquad
\end{equation}
%\vspace{-1mm}
where $Id$ denotes the identity map,  $\Delta$ denotes the Laplacian operator, and $c$ is a positive constant controlling the strength of the preconditioning\footnote{When $c=0$ and $\infty$, the preconditioned gradient is the $L^1$ and $H^1$ gradient, respectively.} (typically $c=1$).
Since the operator ($Id-c\Delta$) is positive definite, the preconditioned gradient (Figure~\ref{fig:grad-prec-arrow} bottom) is still a descent direction. 
While the change of the objective becomes slower w.r.t. flow time\footnote{Here, flow time refers to the time variable $(t)$ in~\eqref{eq:flow-eq}, not the wall-clock time.} (Figure~\ref{fig:grad-precondition} left), much larger time steps are  allowed during line search, thus accelerating  the optimization w.r.t. iteration times (Figure~\ref{fig:grad-precondition} right).    

\begin{figure}[t]
	\centering
	\begin{overpic}[width=0.49\textwidth,keepaspectratio]{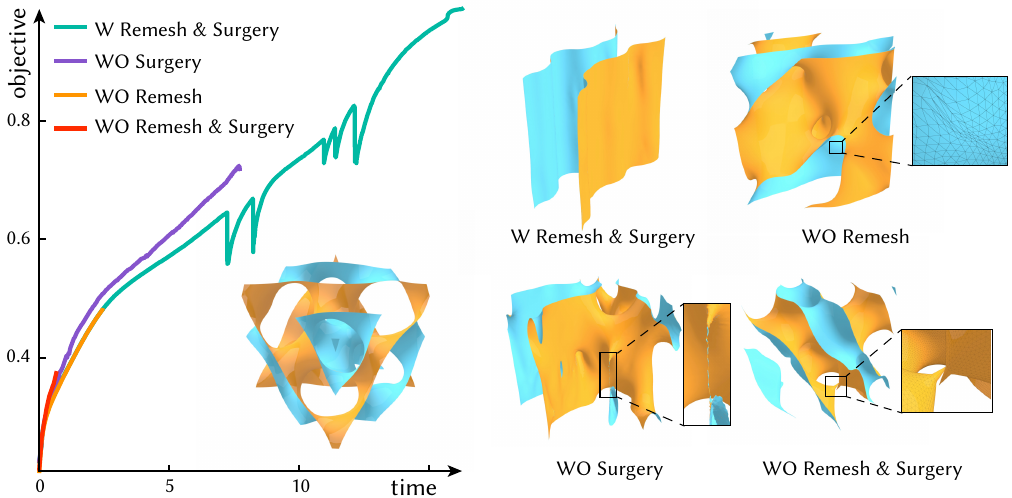}
	\end{overpic}
	%	\vspace{-7mm}
	\caption{
		Maximizing ADS objective $C_A^{3333}$ with/without remeshing and surgery operations.
		The input surface is shown in the left inset.
		The decrease in the objective is observed due to either the surgery operation or failure of the line search (WO surgery).
		Without remeshing or surgery, the optimization stagnates because of poor mesh quality or the formation of necking singularities (right inset).
	}
	%	\vspace{-4mm}
	\label{fig:remesh-willmore}
\end{figure}

\begin{figure}[t]
	% \centering
	\begin{overpic}[width=0.95\linewidth]{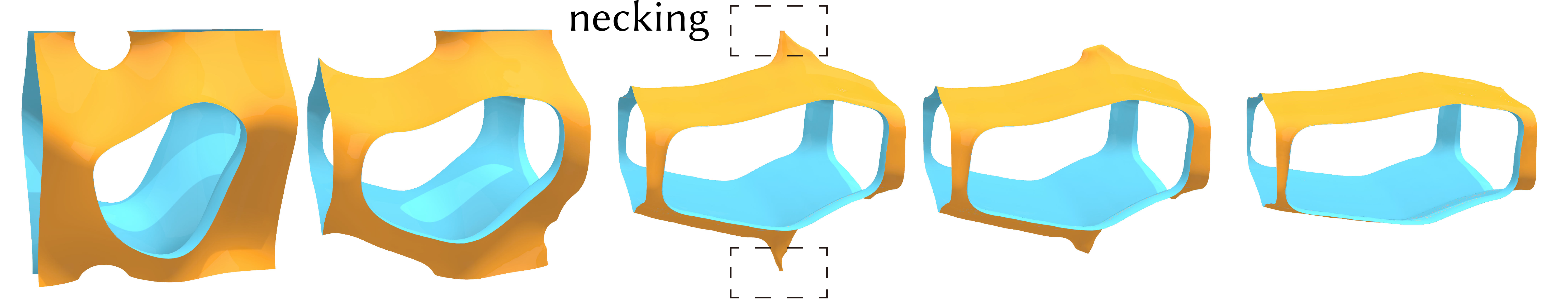}
	\end{overpic}
	%	\vspace{-5mm}
	\caption{
		A necking singularity emerges when maximizing the stiffness under shear strain $\e^{12}$.
		The surface escapes this configuration through numerical surgery.
	}
	%    \vspace{-1mm}
	\label{fig:singularity}
\end{figure}

\begin{figure}[h]
	\centering
	\begin{overpic}[width=0.45\textwidth,keepaspectratio]{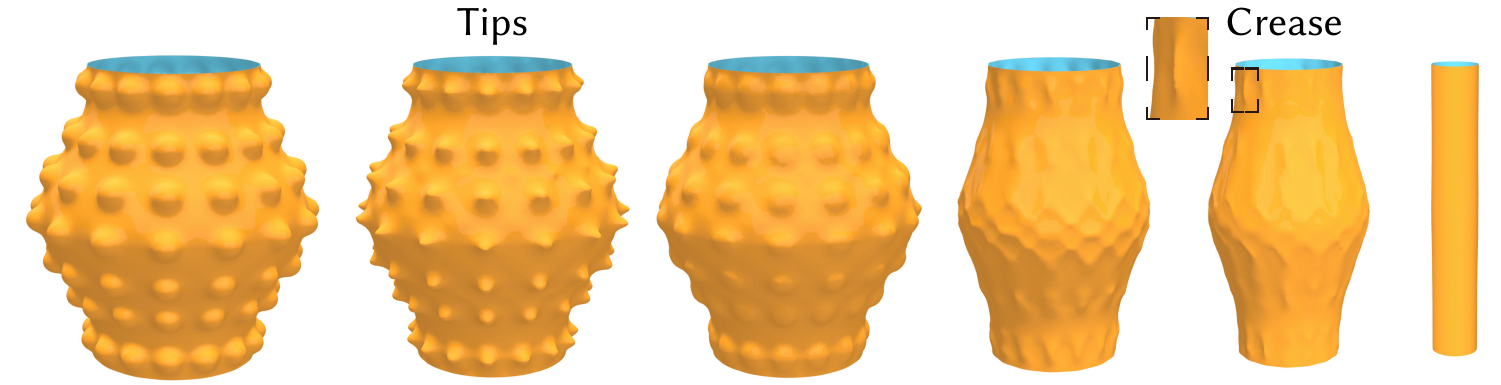}
	\end{overpic}
	%	\vspace{-7mm}
	\caption{
		Tips and creases appear when maximizing Young's modulus along the $z$-axis, and are removed by surgery operations.
		The input and optimized surfaces are shown on the left and right, respectively. 
	}
	%	\vspace{-2mm}
	\label{fig:remove-tip}
\end{figure}

In the discrete setting, we approximate $\grad I_A$ by a vector $V^h$, assembled from the coefficients of per-vertex normal velocities ($v_n^i$) in the discretized form of~\eqref{eq:common-ad-rate} (Appendix~\ref{app:sens}).
We then solve for the preconditioned gradient $\bs d^h$ from the following equation
%\vspace{-1mm}
\begin{equation}
	\label{eq:precond-descent}
	(\mathbf M-c\mathbf L) \bs d^h =-(c+1)V^h,
\end{equation}
%\vspace{-1mm}
where $\mathbf M$ denotes the mass matrix and $\mathbf L$ is the cotangent Laplacian matrix.
Finally, we update the vertex positions by
\begin{equation}
	\bs{x}_i^{k+1} =\bs{x}_i^k + \Delta t \bs{d}^h_i, \quad\forall i\in \mathcal V,
\end{equation}
where $\Delta t$ is a time step, determined using Armijo line search (Alg.\ref{alg:line-search}).
We consider the optimization as converged when either the time step is less than $ 10^{-4}$  for five consecutive iterations, or the time derivative of the objective, approximated by the slope of a linear regression fit to the history of the last 50 iterations, is less than $10^{-3}$.

%\vspace{-2mm}

\begin{algorithm}[t]
	\caption{Armijo line search}
	\label{alg:line-search}
	\LinesNumbered
	\KwIn{Objective $f$, step vector $d_n$, gradient $g_n$}
	\KwOut{Time step $\Delta t$.}
	$\alpha\gets 0.1,\,\tau\gets 0.7$;\\
	$\Delta t \gets AdaptiveTimeStep()$;\\
	\For{$i\gets0$ \KwTo max_iter }{
		\lIf{$f(\mathcal V+\Delta t d_n)<f(\mathcal V)+ \alpha\Delta t d_n^\top g_n$}{Output $\Delta t$}
		$\Delta t\gets \tau\Delta t$\;
	}
	Output $\Delta t$;
\end{algorithm}

\begin{algorithm}[t]
	\caption{Shape optimization for asymptotic properties}
	\label{alg:opti-asym}
	\LinesNumbered
	\KwIn{Periodic mesh $\mathcal M$, objective $f$;}
	\KwOut{Optimized mesh.}
	$i\gets 0$;\\
	$converged\gets false$;\\
	\While{$!converged$ \&\&  i<max_iter}{
		\If{$i\ \%\ 4\ ==\ 0$}{
			$\mathcal M\gets NumericalSurgery(\mathcal M)$;\\
		}
		$\mathcal M\gets DynamicRemesh( \mathcal M)$;\\
		Assemble and solve equation~\eqref{eq:surf-equation-stif};\\
		Evaluate  ADS tensor $\mathbf C_A$
		\tcp*[l]{Section~\ref{sec:the-vari-eq} }
		Evaluate objective $f^i$;\\
		Evaluate the $L^2$ gradient $g_n$ based on~\eqref{eq:custom-obj-der},~\eqref{eq:l2-graident},~and~\eqref{eq:ca-ijkl-dt};\\
		Solve for preconditioned gradient $d_n$ from~\eqref{eq:precond-descent};\\
		$\Delta t \gets ArmijoLineSearch(f,\ g_n,\ d_n)$
		\tcp*[l]{Alg.\ref{alg:line-search}}
		$\mathcal V\gets \mathcal V+ \Delta t d_n $;\\
		$converged\gets CheckConvergence(f^i, \Delta t)$;\\
		$i\gets i+1$;\\
	}
	Output $\mathcal{M}$;
\end{algorithm}

\paragraph{Adaptive initial time step}
\label{sec:adaptive-set}
During the line search, each evaluation of objective requires solving equation~\eqref{eq:surf-equation-stif}, which is time consuming.
A large initial time step leads to excessive backtracking steps and increased  equation solves.
To reduce the overhead, we set the initial step in the Armijo line search as the average of those determined in the last five iterations, multiplied by a factor of 2 to avoid reducing the time step to zero.

%\vspace{-2mm}

\paragraph{Dynamic remesh}
\label{sec:dyn-remesh}
As the mesh quality degrades during the flow, the optimization eventually stalls due to the increasingly small time step, even with the preconditioned gradient.
To address this, we  employ the dynamic remeshing algorithm proposed in~\cite{Yu2021-srf}.
We refer readers to Section 6 in~\cite{Yu2021-srf} for implementation details.
With the remeshing procedure in place, the optimization proceeds normally, see the comparison in Figure~\ref{fig:remesh-willmore}.

%\vspace{-2mm}
\paragraph{Numerical surgery}
\label{sec:surg}
As the optimization progresses, the surface may develop necking singularities  (Figure~\ref{fig:singularity}).
Inspired by~\cite{kovacs2024-surgery}, we apply a numerical surgery procedure to remove them.
We identify  candidate regions based on the maximal principal curvature, the triangles with vertex curvature values larger than a threshold  (typically 25) are selected and grouped into  connected components.
For each component, we first remove its constituent triangles and fill the resulting hole using the CGAL~\cite{cgal:h-af-24b} library, then  a bi-Laplacian fairing step~\cite{Botsch04-fairing} is performed in the $k$-ring of the filling patch, where $k=4$ in our experiments.
This procedure also eliminates undesired tips and creases that may arise during the optimization (Figure~\ref{fig:remove-tip}).
See also the ablation study in Figure~\ref{fig:remesh-willmore}.
To reduce the computational overhead, this operation is performed every four iterations.

%\vspace{-2mm}
\subsubsection{Custom objective}
With above techniques,  we can optimize a custom objective function $f(\mathbf C_A)$, with sensitivity given by
\begin{equation}
	\label{eq:custom-obj-der}
	\frac{d}{dt}f = \frac{\partial f}{\partial C_{A}^{ijkl}}\dt{C}_A^{ijkl},
\end{equation}
where we adopted the convention of Einstein summation.
The  gradient evaluation and optimization steps remain the same as before.
Pseudocodes of our overall optimization pipeline are shown in Alg.~\ref{alg:opti-asym}.

%	\vspace{-5mm}

\section{Validation of theory}
\label{sec:validation}
In this section, we first illustrate our method for generating and handling periodic surfaces (Section~\ref{sec:handle-period-surf} and~\ref{sec:gen-period-surf}).
Then we conduct numerical experiments in Section~\ref{sec:num-problem} to verify the convergence of ADS formula~\eqref{eq:dir-asymp-stif} and assess the accuracy of our discretization scheme.
In Section~\ref{sec:validate-tpms-opt}, we validate the optimality of TPMS as concluded in Section~\ref{sec:opti-tpms-bulk}.

\begin{figure}[t]
	\begin{overpic}[width=0.90\linewidth]{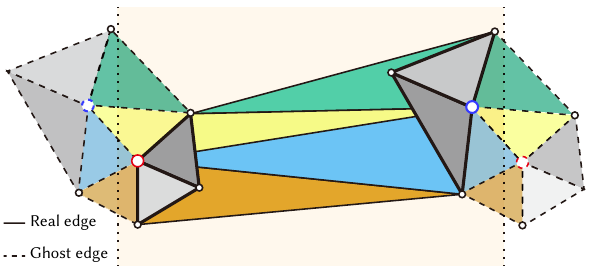}
		\put(16,46){$x=-1$}
		\put(81,46){$x=1$}
		\put(22,21) {$v_i$}
		\put(87,21) {$v_i$}
		\put(18,29) {$v_j$}
		\put(82,29) {$v_j$}
		\put(55,21) {$e_{ij}$}
	\end{overpic}
%	\vspace{-4mm}
	\caption{
		Reproduced 1-ring of vertex ($v_i$ and $v_j$) when the incident edge is longer than a threshold.
	}
	\label{fig:period-mesh}
\end{figure}
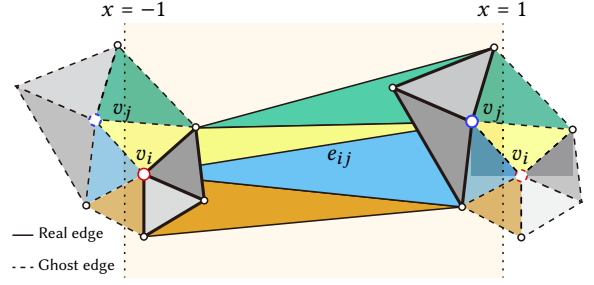

\subsection{Handling periodic mesh}
\label{sec:handle-period-surf}

To facilitate the query of topological information on the mesh, such as  incident vertices or faces, we merge each pair of periodic vertices into a single vertex. 
The position of this vertex is taken from either of the two original vertices.
This merging operation may result in elongated edges, which we exploit to identify whether an incident vertex belongs to another period (Figure~\ref{fig:period-mesh}). 
Specifically, when accessing the incident vertices of a given vertex, we check the difference vector between their positions.
If any component of the vector exceeds 1 or falls below -1, we translate the incident vertex along that axis by $2k$ or $-2k$ until the difference lies within the range $[-1,1]$ (Figure~\ref{fig:period-mesh}).
These translated vertices are only used on the fly and are not stored back to the mesh data structure.
With such operation, the mesh becomes locally equivalent to the periodic one. 
When updating vertex positions during the flow (Section~\ref{sec:geo-opt}), only a small translation vector is applied to each vertex, which does not interfere with the aforementioned operation. 
While long edges not caused by merging could potentially pose a problem,
we avoid this by keeping the edge lengths less than $0.3$ during the remeshing process (Section~\ref{sec:dyn-remesh}).
Note that the remeshing algorithm involves only local operations, such as vertex translation, edge collapse (split), and edge flipping, this operation is not affected either.

\subsection{Generation of periodic surface}
\label{sec:gen-period-surf}
\begin{figure}[t]
	\centering
	\begin{overpic}[width=0.45\textwidth,keepaspectratio]{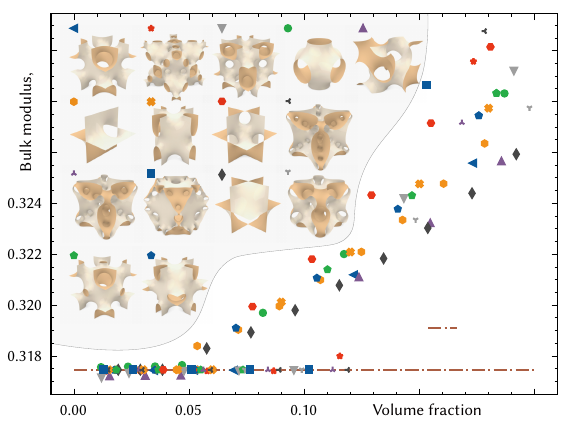}
		\put(84,2) {\small$\rho_\eps$}
		\put(3,62) {\rotatebox{90}{\small$K_\eps/\rho_\eps$}}
		\put(81,16) {$K_A^*$}
	\end{overpic}
	%	\vspace{-3mm}
	\caption{
		Scatter plot of the the normalized bulk modulus ($K_\eps/\rho_\eps$) of 15 types of TPMS (inset) as the thickness approaches zero. 
		The upper bound $K_A^*\approx 0.3174$ in our setting.
	}
	\label{fig:tpms-scat}
\end{figure}

To validate the conclusions in Section~\ref{sec:opti-tpms-bulk}, we generate 15 types of TPMS (Figure~\ref{fig:tpms-scat} inset) using Surface Evolver~\cite{surf-evolver}, and sample 1000 non-TPMS surfaces for numerical experiments.
All surfaces are rescaled to $Y=[-1,1]^3$.
The non-TPMS surfaces are sampled by perturbing the trigonometric approximations of Schwarz P, Schwarz D, IWP, and Gyroid, which are the 
level sets of the following functions:
\begin{equation}
	\label{eq:trig-approx}
	\small
	\begin{aligned}
		&\varphi_P(\boldsymbol{r})=  \cos (2 \pi x)+\cos (2 \pi y)+\cos (2 \pi z)\\
		&\varphi_G(\boldsymbol{r})=  \sin (2 \pi x) \cos (2 \pi y)+\sin (2 \pi z) \cos (2 \pi x)+\sin (2 \pi y) \cos (2 \pi z)\\
		&\varphi_{D}(\boldsymbol{r})=  \sin (2 \pi x) \sin (2 \pi y) \sin (2 \pi z) + \sin (2 \pi x) \cos (2 \pi y) \cos (2 \pi z)\\&+ \cos(2\pi x)\sin(2\pi y)\cos(2\pi z)+\cos(2\pi x)\cos(2\pi y)\sin(2\pi z)\\
		&\varphi_{\mathrm{IWP}}(\boldsymbol{r}) = 2[\cos (2 \pi x) \cos (2 \pi y)+\cos (2 \pi y) \cos (2 \pi z)\\&+\cos (2 \pi z) \cos (2 \pi x)]-[\cos (2 \cdot 2 \pi x)  +\cos (2 \cdot 2 \pi y)+\cos (2 \cdot 2 \pi z)].
	\end{aligned}
\end{equation}
We perturb their level sets by adding an additional noise term:
$
	\delta\varphi = \sum_i s_i b^i(\bs x),
$
where $b^i(\bs x)$ denotes the basis functions; $s_i$ are the coefficients randomly sampled from a uniform distribution over $[-\bar s, \bar s]$, with the parameter $\bar s$ used to control the perturbation strength.
The basis functions are drawn from the set 
$B_N:=T_N \cup \br{pq: p\in T_N, q\in T_N}$, 
where 
$T_N =\br{\sin(2\pi k x_i),\,\cos(2\pi k x_i)}_{k=1,\cdots,N}$.
Here, integer $N$ controls the degrees of freedom (DoF) in the perturbation.
We observe that choosing $N \ge 2$ often produces overly noisy surfaces; therefore, we set $N = 1$ in our experiments.
We evaluate the perturbed function 
$
	\tilde{\varphi} := \varphi + \delta \varphi
$
%\vspace{-1mm}
on a  $128\times128\times128$ grid evenly discretizing domain $[0,1]^3$, where $\varphi$ denotes  any of the functions in~\eqref{eq:trig-approx}.
The level set surface is extracted from the grid values using the Marching Cube algorithm~\citep{marching-cube} implemented in \verb|libigl|~\cite{libigl}.
We then improve the mesh quality using the incremental isotropic remeshing algorithm~\cite{pmp-book}.
We set the perturbation strength $\bar s =0.1$, $0.2$, $0.3$, $0.4$, $0.5$, and for each setting, we sample the coefficients $\br{s_i}$ for $50$ times, generating 50 perturbed surfaces.
These surfaces are used to validate the ABM upper bound (Section~\ref{sec:validate-tpms-opt}) and as inputs to our optimization algorithm.

\subsection{Simulation setup}
We import the generated surfaces to ABAQUS for finite element simulation. 
Each triangular facet of the surface mesh is assigned a STRI3 shell element for mechanical simulation.
The Young's modulus and Poisson's ratio are set to $Y=1 $ and $\nu =0.3$, respectively. 
We use the Micromechanics plugin to enforce periodic boundary conditions and compute the effective elastic tensor.

\subsection{Convergence test}
\label{sec:num-problem}

\begin{figure}[t]
	\centering
	\begin{overpic}[width=0.45\textwidth,keepaspectratio]{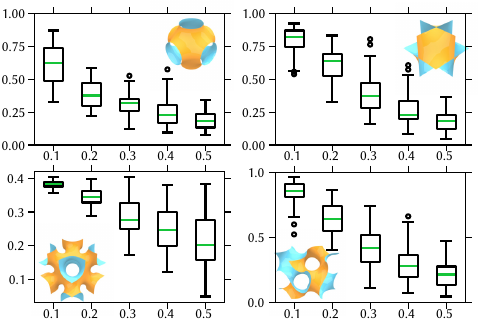}
		\put(-3,45) {\rotatebox{90}{\tiny$K_A/ K_A^*$}}
		\put(-3,14) {\rotatebox{90}{\tiny$K_A/ K_A^*$}}
		\put(26,-3) {$\bar s$}
		\put(77,-3) {$\bar s$}
	\end{overpic}
	\caption{
		The ABM of perturbed surfaces derived from Schwarz P, Schwarz D,
		IWP, and Gyroid. Each subplot corresponds to one TPMS type. The abscissa
		denotes the perturbation strength $\bar s$, with 50 surfaces sampled for each
		setting.
		Here, $K_A^*$ denotes the upper bound in Theorem~\ref{thm:tpms-opt-bulk}.
	}
	%	\vspace{-3mm}
	\label{fig:sample-surf-bulk}
\end{figure}
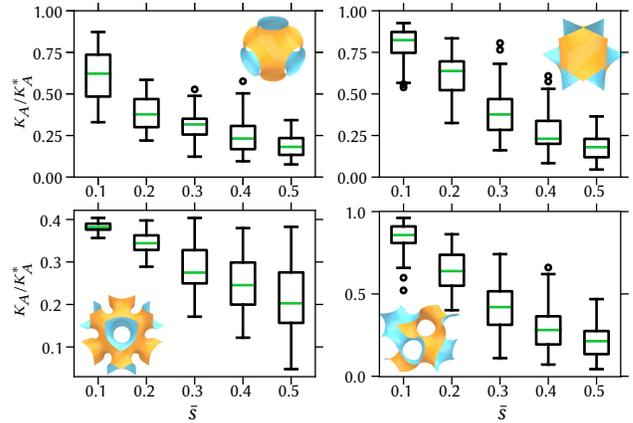
According to the definition of asymptotic elastic tensor~\eqref{eq:EA-tensor-def}, it suffices to validate the following relation
\begin{align}
	\label{eq:valid-CA-kA}
	\mathbf C_A=\lim_{\eps\to0}\frac{{\mathbf C}_\eps}{\rho_\eps(\tw)},
\end{align}
where $\mathbf C_A$ is defined in~\eqref{eq:ca-def}.
To assess the convergence,  we introduce the following relative error
\begin{align}
	\label{eq:err-ads}
	\mathcal{E}^s &:= \frac{\|\mathbf C_A^h-\mathbf C_\eps/\rho_\eps\|_F}{\|\mathbf C_A^h\|_F},
\end{align}
where $\mathbf C_A^h$ is given in~\eqref{eq:disc-ca-form}, and $\mathbf C_\eps$ is obtained from the simulations in ABAQUS.
We generate the following level set surface for test    
\begin{equation}
	\varphi(\bs r) = 1.2\cos(\pi x) + 0.9\cos(\pi y) + 1.5\cos(\pi z).
\end{equation}
\begin{wrapfigure}{r}{0.42\columnwidth}
	\hspace{-10pt}
	\vspace{-14pt}
	% \centering
	\begin{overpic}[width=0.99\linewidth]{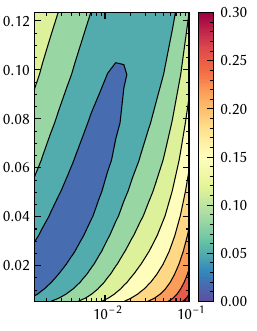}
		\put(18,3){$2\eps$}
		\put(5,98){$h$}
		\put(60,100){$\mathcal E^s$}
	\end{overpic} 
	\caption{
	}
	\vspace{-6mm}
	\label{fig:converg-mech}
\end{wrapfigure}
The surface mesh is extracted using the method in Section~\ref{sec:gen-period-surf}.
We refine it using the remeshing algorithm in Section~\ref{sec:dyn-remesh},  setting the target edge length to $0.1\times (4/5)^k$, where $k$ denotes the refined times.
After each refinement, the vertices are projected to the level set through a binary search along the gradient direction.
We then evaluate  $\mathbf C_A^h$ and the corresponding average circumradius of the mesh faces (considered as element size $h$).
On the other hand, we import mesh to ABAQUS for simulation, with a  decreasing shell thickness $2\eps = 0.1\times(4/5)^k$. 
The error $\mathcal E^s$ for each pair of $h$ and $2\eps$ is plotted in Figure~\ref{fig:converg-mech}.
\begin{wrapfigure}{r}{0.42\columnwidth}
	\hspace{-10pt}
	\vspace{-8pt}
	\begin{overpic}[width=0.99\linewidth]{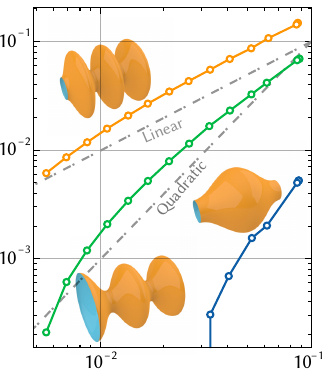}
		\put(55,0){\small $h$}
		\put(5,100){\small $\mathcal{E}^h$}
	\end{overpic} 
	\vspace{-1mm}
	\caption{
	}
	\vspace{-3mm}
	\label{fig:revo-conv}
\end{wrapfigure}
It shows that the error tends to zero when $h$ and $\eps$ decrease simultaneously. 

Next, we construct three revolution surfaces around the $x$-axis (Figure ~\ref{fig:revo-conv} inset), for which the ADS under strain $\e^{11}$ can be solved semi-analytically (SM Section~5).
Then we compare the ADS ($E_A^h$) obtained from our approach with the semi-analytical solution  ($E_A$) under strain $\e^{11}$.
We evaluate the relative error 
%\vspace{-1mm}
\begin{equation}
	\mathcal{E}^h := \frac{|E_A^h-E_A|}{|E_A|}\qquad\qquad\qquad \qquad\qquad\qquad
\end{equation}
using a decreasing element size~($h$).
The experimental results (Figure ~\ref{fig:revo-conv}) show a linear or even better convergence rate.

While we have no theoretical guarantee for the convergence of this discretization scheme,  we argue that the variational problem~\eqref{eq:surf-vari-form} determining ADS is fundamentally difficult to discretize with guaranteed convergence.
The reason is twofold.
First, the function space is generally too large, possibly beyond the space of distributions, as noted in~\cite{CiarletGenmem1996,meunier2006}.
This leaves the linear form (RHS of~\eqref{eq:surf-vari-form}) in a very small dual space $(\dot{V}_\#(\tw))'$, making it unclear whether the discretized form still defines a valid functional in this space, or whether the discretization can approximate the solutions that may not be regular functions.
Second, the ADS could be extremely sensitive to the geometric perturbations in special cases, where an infinitesimally small error in the geometry could result in a significant change in the ADS (SM Section~6).

\subsection{The optimality of TPMS}
\label{sec:validate-tpms-opt}

We validate the conclusions in Section~\ref{thm:tpms-opt-bulk}.
For the generated 15 types of TPMS (Section~\ref{sec:gen-period-surf}), we approximate $\mathbf C_A$ by $\mathbf C_\eps/\rho_\eps$ with a decreasing thickness, 
and use identity $K_A=\frac{1}{9}\mathbf C_A:\mathbf I:\mathbf I$ to compute ABM.
The approximated ABMs and their corresponding volume fractions are shown in Figure~\ref{fig:tpms-scat}.
The results indicate that the ABMs are very close to the upper bound, even when the volume fraction is relatively large ($<20\%$).

On the other hand, the ABMs of non-TPMS surfaces are approximated in the same way with a fixed thickness of $0.004$.
The simulation results for the relative magnitudes $K_A/K_A^*$ are shown in Figure~\ref{fig:sample-surf-bulk}.
We observe that $K_A/K_A^*\le 1$ for all these surfaces.
Moreover, as the perturbation strength increases, this  value tends to decrease as they are getting distinct from TPMS\footnote{
	The ABMs of the perturbed IWP surfaces are exceptionally low, which may be attributed to either the poor accuracy of the trigonometric representation of IWP or the sensitive problem of the ADS, as discussed in Section~\ref{sec:num-problem}.
}.
These numerical results conform to our theoretical prediction.

\section{Validation of the optimization algorithm}

\subsection{Optimize ADS}
\label{sec:tailor-ads}

\begin{figure}[h]
	\centering
	\begin{overpic}[width=0.43\textwidth,keepaspectratio]{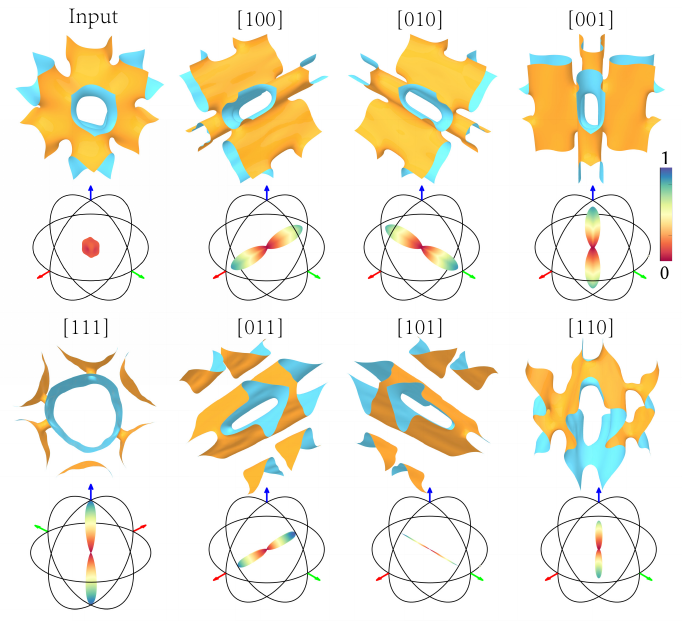}
	\end{overpic}
%	\vspace{-3mm}
	\caption{
		Distribution of Young's modulus before and after optimization, visualized by VELAS~\cite{velas2023}.
		The objective directions ($\mathbf{z}$) of Young's modulus are noted above the surfaces.
	}
%	\vspace{-1mm}
	\label{fig:ads-max-G}
\end{figure}

\paragraph{Young's modulus}
We maximize Young's modulus in different directions with a single input surface.
The directional Young's modulus~\cite{Schumacher2018} is given by
\begin{equation}
	\label{eq:ads-young-z}
	Y_{\mathbf z}(\mathbf C_A)=  \pr{\bs{ \sigma}_{\mathbf z}:(\mathbf C_A)^{-1}:\bs \sigma_{\mathbf z}}^{-1},
\end{equation}
where $\bs\sigma_{\mathbf z}=\mathbf{z}\mathbf{z}^\top$ and $\mathbf z$ denotes the direction in which Young's modulus is computed. 
We use  \verb|AutoDiff|  module of the Eigen library~\cite{eigenweb} to compute its derivative.
The Young's modulus distributions of the surfaces before and after optimization are shown in Figure~\ref{fig:ads-max-G}.
A significant increase in stiffness is observed for various objective directions after optimization.

\begin{figure}[h]
	\centering
	\begin{overpic}[width=0.43\textwidth,keepaspectratio]{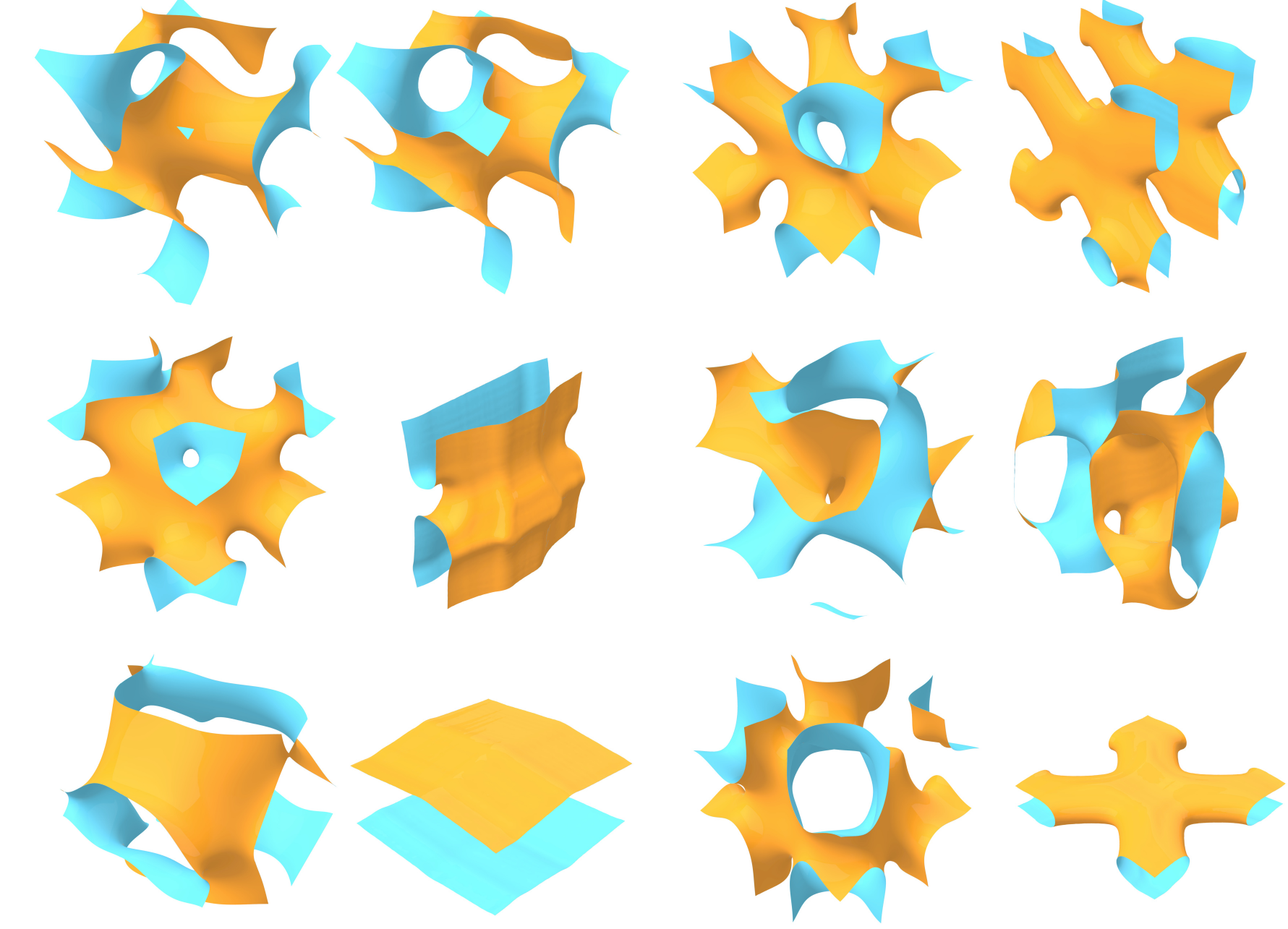}
		\put(10,73){\small Input}
		\put(31,73){\small Result}
		\put(60,73){\small Input}
		\put(82,73){\small Result}
		% % 
		\put(25,49) {(a)}
		\put(75,49) {(b)}
		\put(25,22) {(c)}
		\put(75,22) {(d)}
		\put(25,-2) {(e)}
		\put(75,-2) {(f)}
		\put(-1,47) {\footnotesize $C^{2323}_A$}
		\put(10,47) {\footnotesize $0.0950$}
		\put(33,47) {\footnotesize $0.1932$}
		\put(60,47) {\footnotesize $0.0683$}
		\put(85,47) {\footnotesize $0.2699$}
		\put(-1,22) {\footnotesize $C^{1313}_A$}
		\put(10,22) {\footnotesize $0.0637$}
		\put(33,22) {\footnotesize $0.3114$}
		\put(60,22) {\footnotesize $0.0457$}
		\put(85,22) {\footnotesize $0.1988$}
		\put(-1,-2) {\footnotesize $C^{1212}_A$}
		\put(10,-2) {\footnotesize $0.0428$}
		\put(33,-2) {\footnotesize $0.3723$}
		\put(60,-2) {\footnotesize $0.0546$}
		\put(85,-2) {\footnotesize $0.2806$}
	\end{overpic}
%	\vspace{-1mm}
	\caption{
		Maximizing stiffness under three shear strains $\e^{23}$ (Top), $\e^{13}$ (Middle), and $\e^{12}$ (Bottom). 
		The corresponding values of the objectives are noted below the surfaces.
	}
	\label{fig:shear-opt}
\end{figure}
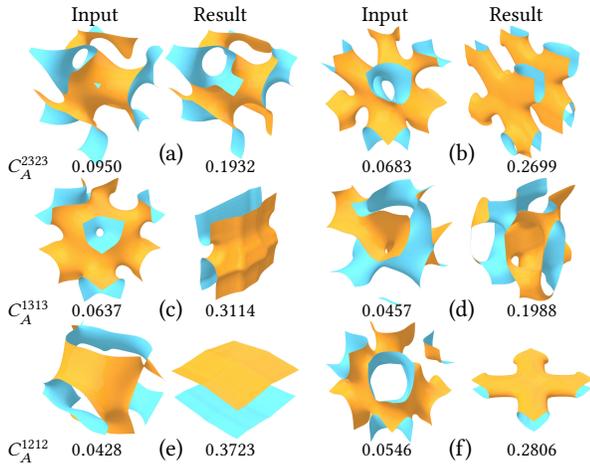

\begin{figure}[t]
	\centering
	\begin{overpic}[width=0.99\linewidth]{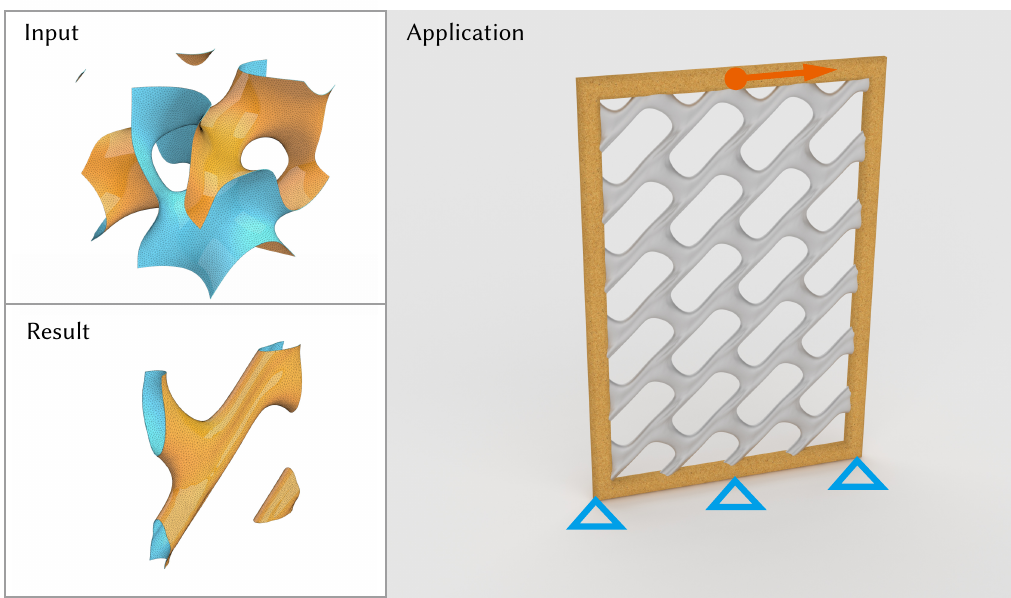}
	\end{overpic}
%	\vspace{-3mm}
	\caption{
		Steel tube with optimized shear modulus is used to strengthen a frame.
		The red arrow denotes the load.
	}
%	\vspace{-1mm}
	\label{fig:app-fence}
\end{figure}

\paragraph{Shear modulus}
Next, we maximize the stiffness under shear strains $\e^{23}$, $\e^{13}$, and $\e^{12}$ with different input surfaces.
Thanks to the surgery operation, the optimization removes the topological barriers  and yields a better result (Figure~\ref{fig:shear-opt}).
See also an application in Figure~\ref{fig:app-fence}.

\begin{figure}[h]
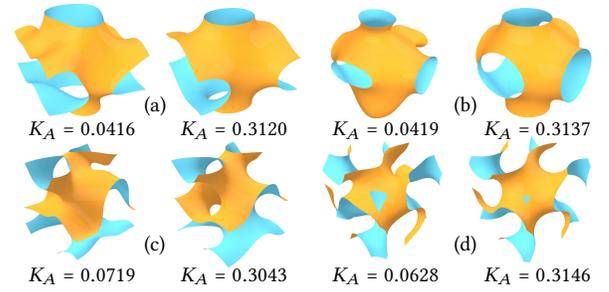

%	\vspace{-1mm}
	\centering
	\begin{overpic}[width=0.45\textwidth,keepaspectratio]{bulk-tpms1-c.pdf}
		\put(5, -2){\small $K_A=0.0719$}
		\put(30,-2){\small $K_A=0.3043$}
		\put(55,-2){\small $K_A=0.0628$}
		\put(80,-2){\small $K_A=0.3146$}
		\put(5, 22){\small $K_A=0.0416$}
		\put(30,22){\small $K_A=0.3120$}
		\put(55,22){\small $K_A=0.0419$}
		\put(80,22){\small $K_A=0.3137$}
		\put(24, 26){\small (a)}
		\put(75,26){\small (b)}
		\put(24,3){\small (c)}
		\put(75,3){\small (d)}
	\end{overpic}
%	\vspace{-2mm}
	\caption{
		Surfaces converge to TPMS-like shapes when optimizing ABM.
		The first and third columns show the input surfaces, and second and fourth columns show the optimized results.
	}
	\label{fig:abm-opt}
\end{figure}

\paragraph{Bulk modulus}
As Theorem~\ref{thm:tpms-opt-bulk} indicates, the upper bound of ABM is $K_A^*=\frac{4}{9}\pr{\lambda_0+\mu}\approx 0.3174$ under our experimental setting, and is achieved when the surface is a TPMS.
Therefore, we conduct experiments with ABM as the objective to observe how the surfaces deform after optimization.
The results are shown in Figure~\ref{fig:abm-opt}.
Interestingly, the surfaces converge toward TPMS-like shapes with an increased ABM,  consistent with our theoretical predictions.

\begin{figure}[h]
	\centering
%	\vspace{2mm}
	\begin{overpic}[width=0.45\textwidth,keepaspectratio]{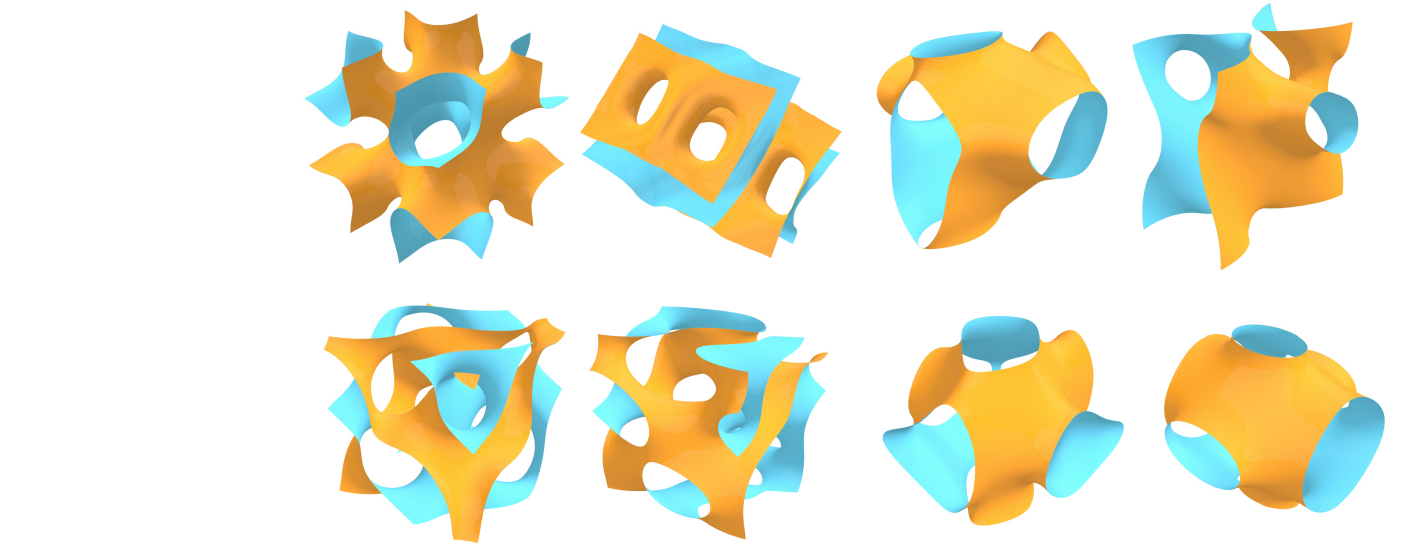}
		\put(-2,27){\small $\begin{pmatrix}
				1 & 0 & 1\\
				0 & -2& 0\\
				1 & 0 & 1
			\end{pmatrix}$ }
		\put(-1,10){\small $\begin{pmatrix}
				1 & 0 & 0\\
				0 & 2& 0\\
				0 & 0 & 1
			\end{pmatrix}$ }
		\put(7,38) {$\veps$}
		\put(27,39) {\small Input}
		\put(44,39) {\small Result}
		\put(65,39) {\small Input}
		\put(82,39) {\small Result}
		\put(22,18.5) {\footnotesize $E_A=0.7574$}
		\put(42,18.5) {\footnotesize $E_A=4.4476$}
		\put(63,18.5) {\footnotesize $E_A=0.5285$}
		\put(83,18.5) {\footnotesize $E_A=3.5150$}
		\put(39,22) {\footnotesize (a)}
		\put(78,22) {\footnotesize (b)}
		\put(39,2) {\footnotesize (c)}
		\put(78,2) {\footnotesize (d)}
		\put(22,-2) {\footnotesize $E_A=1.8198$}
		\put(42,-2) {\footnotesize $E_A=5.5327$}
		\put(63,-2) {\footnotesize $E_A=1.5003$}
		\put(83,-2) {\footnotesize $E_A=6.2003$}
	\end{overpic}
	\caption{
		Optimization of stiffness under specific strains (noted on the left side).
	}
%	\vspace{-3mm}
	\label{fig:vari-strain}
\end{figure}
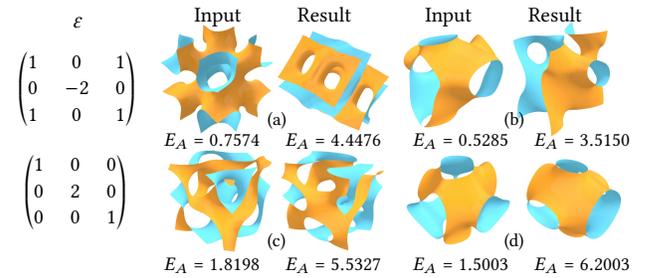

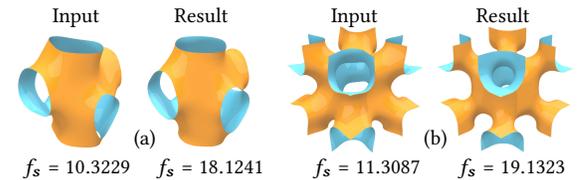
\begin{figure}[h]
	\centering
%	\vspace{2mm}
	\begin{overpic}[width=0.43\textwidth,keepaspectratio]{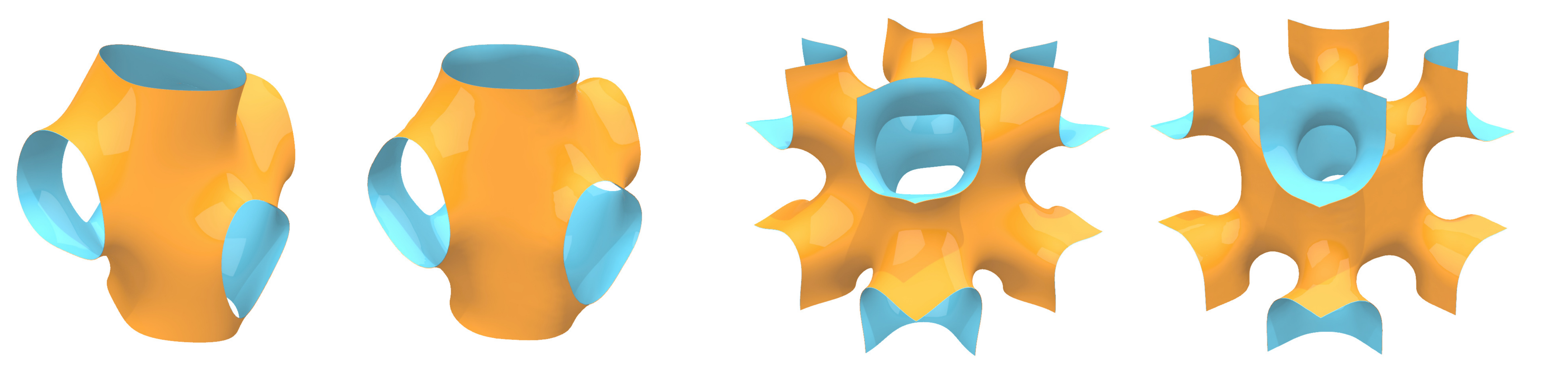}
		\put(7,24){\small Input}
		\put(28,24){\small Result}
		\put(55,24){\small Input}
		\put(80,24){\small Result}
		\put(2,-2){\small $f_{\bs s}=10.3229$}
		\put(25,-2){\small$f_{\bs s}=18.1241$}
		\put(52,-2){\small$f_{\bs s}=11.3087$}
		\put(77,-2){\small$f_{\bs s}=19.1323$}
		\put(21,3){\small (a)}
		\put(71,3){\small (b)}
	\end{overpic}
	\caption{
		Optimizing ADS under the deviatoric strain $\bs s = \operatorname{diag}(1, -2, 1)$ over all spatial orientations. Each result shows an increase in the objective value $f_{\bs s}$ after optimization.
	}
%	\vspace{-3mm}
	\label{fig:dev-strain}
\end{figure}

\paragraph{Arbitrary strain}
We can apply our approach to maximize the stiffness under arbitrary strains.
Figure~\ref{fig:vari-strain} shows the improvement in ADS under two representative strains.
Optimizing ADS over a range of strains is also supported.
	For instance, we consider the average ADS under deviatoric strain $\bs s = \operatorname{diag}(1, -2, 1)$, integrated over all spatial orientations.
	The corresponding objective is:
	\begin{equation}
		f_{\bs s}(\mathbf C_A)=\int_{SO(3)}\veps(\mathbf R):\mathbf C_A:\veps (\mathbf R)\ d\mu(\mathbf R)
	\end{equation}
	where $\veps(\mathbf R) = \mathbf R \bs s  \mathbf R^\top$ denotes the rotated strain tensor under the orientation $\mathbf R \in SO(3)$, and $d\mu(\mathbf R)$ is the Haar measure on $SO(3)$.
	The computation (SM Section 7) shows that this objective simplifies to
	\begin{equation}
		\label{eq:dev-obj-test}
		f_{\bs s}(\mathbf C_A)=\sum_{i}C_A^{iiii}+3\sum_{i\ne j} C_A^{ijij}-\sum_{i\ne j}C_A^{iijj}.
	\end{equation}
	A demonstration is shown in Figure~\ref{fig:dev-strain}.

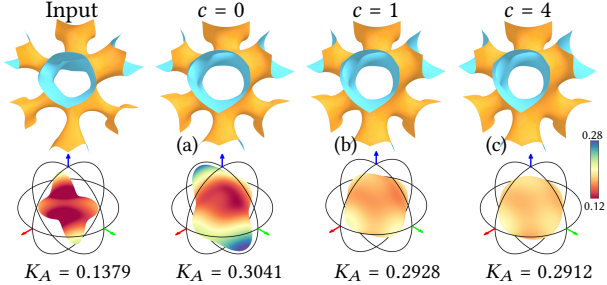
\begin{figure}[h]
	\centering
%	\vspace{1mm}
	\begin{overpic}[width=0.45\textwidth,keepaspectratio]{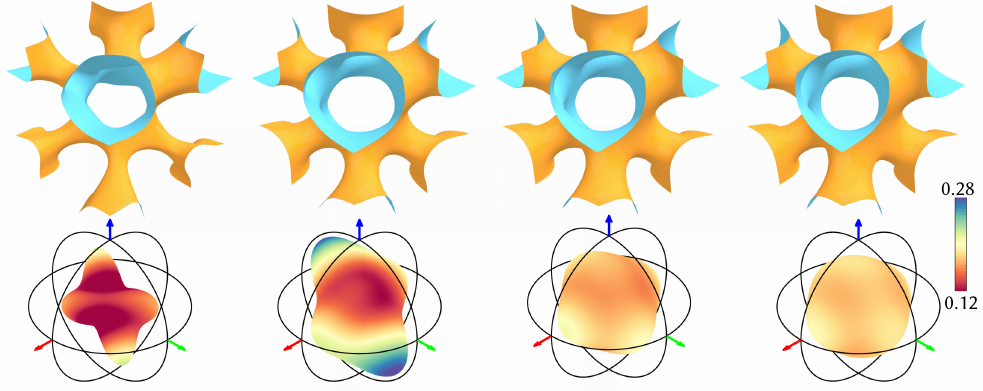}
		\put(7.5, 40){Input}
		\put(32.5,40){$c=0$}
		\put(58,40){$c=1$}
		\put(83,40){$c=4$}
		\put(4, -3){\small$K_A = 0.1379$}
		\put(29,-3){\small$K_A=0.3041$}
		\put(55,-3){\small$K_A= 0.2928$}
		\put(80,-3){\small$K_A=0.2912$}
		\put(29,18){\small(a)}
		\put(55,18){\small(b)}
		\put(80,18){\small(c)}
	\end{overpic}
%	\vspace{-0mm}
	\caption{
		Distribution of Young's modulus after optimization with different isotropic penalty coefficients (top), visualized by VELAS~\cite{velas2023}.
		The corresponding ABM values are shown below.
	}
%	\vspace{-0mm}
	\label{fig:cusobj}
\end{figure}

\paragraph{Custom objective}
Custom objective is also supported.
For instance, we can add an isotropic penalty $f^s_{iso}(\mathbf C_A)$  (Appendix~\ref{sec:iso-penal}) to the objective when optimizing ABM, i.e., we employ the following objective  
\begin{equation}
	f(\mathbf C_A)= \frac{1}{9}\mathbf{I}:\mathbf C_A:\mathbf{I} - c f^s_{iso}(\mathbf C_A),
\end{equation}
where $c$ is a constant penalty coefficient.
A comparison across different values of $c$ is shown in Figure~\ref{fig:cusobj}.

We clarify that our method is not suitable for designing surfaces with a specified asymptotic elastic tensor or for minimizing stiffness under some strains.
The reason is that when the objective encourages more compliant behavior—such as minimizing stiffness or targeting small ADS values—the surface tends to form wrinkles that help reduce stiffness. 
However, this often compromises the smoothness of the mesh, making our approximation inaccurate.
See the example of designing metamaterials with a negative Poisson's ratio\footnote{A negative Poisson's ratio is typically associated with a  compliant material.}  in Figure~\ref{fig:adsnpr}.
The employed objective is elaborated in Appendix~\ref{sec:opt-npr}.

	\begin{remark}
		Although it has been demonstrated that high-frequency surface details can   enhance structural stability~\cite{Montes2023}, such details may introduce manufacturing defects in metamaterials, as the physical size of the unit cell is typically small, making the intricate features difficult to reproduce accurately. Besides, they increase the risk of stress concentration.
			Therefore, despite this being a limitation of ADS, we believe that smoother surfaces are generally more desirable in practice.
	\end{remark}
	\begin{remark}
		Additionally, we believe the shell lattices can still exhibit strong  structural stability, such as buckling resistance, even with a fairly smooth middle surface and as thickness approaches zero, provided that the inextensional displacement space (i.e., \(\ker\gamma\)) degenerates to zero.
			This space contains the most likely buckling deformations, as they induce no membrane strain.
			A possible approach to enhancing structural stability based on this insight is left in the future work, as we will discuss in Section~\ref{sec:conc}.
	\end{remark}

\begin{figure}[h]
	\centering
	\begin{overpic}[width=0.45\textwidth,keepaspectratio]{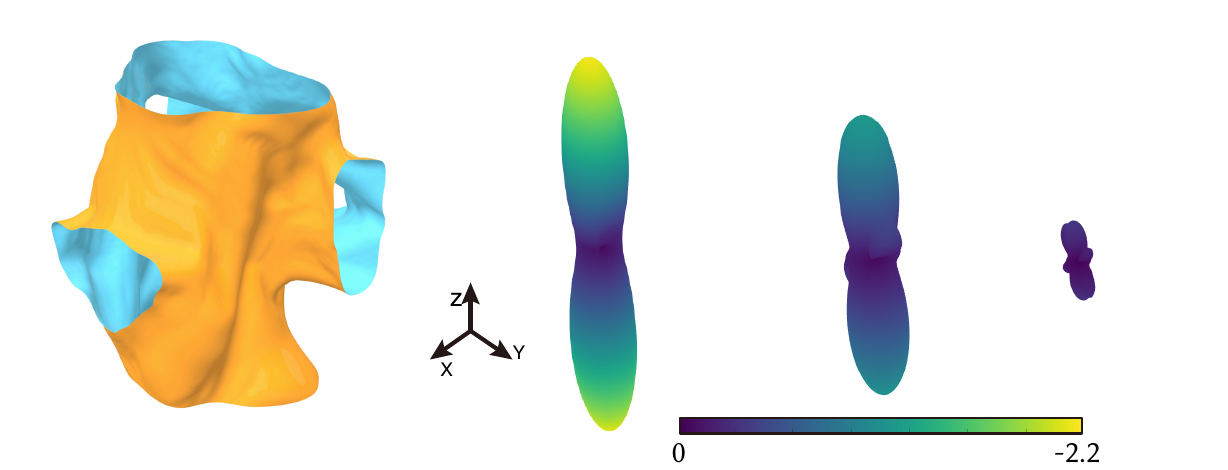}
		\put(43,34){$\nu(\mathbf C_A)$}
		\put(65,34){$\eps=0.02$}
		\put(85,34){$\eps=0.05$}
	\end{overpic}
%	\vspace{-1mm}
	\caption{
		Minimizing the Poisson's ratio in the $z$-direction.
		The wrinkles appear on the optimized surface (Left).  
		The distribution of the negative Poisson's ratio is shown on the right, visualized by VELAS~\cite{velas2023}.
		Although the Poisson's ratio $\nu(\mathbf C_A)$ of the surface attains $-2.2$ (Right first),  it decreases rapidly as the thickness increases (Right second and third). 
	}
%	\vspace{0mm}
	\label{fig:adsnpr}
\end{figure}

\subsection{Evolution of effective properties of shell lattice}
\begin{figure}[t]
	\centering
	\begin{overpic}[width=0.49\textwidth,keepaspectratio]{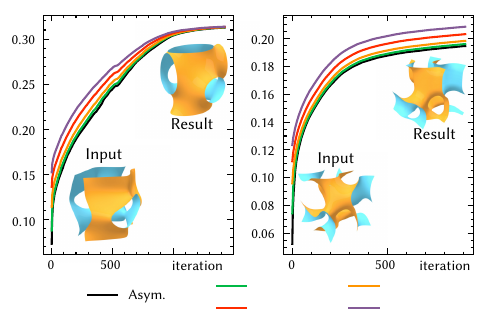}
		\put(52,4.5){\small $\eps =0.01$}
		\put(78,4.5){\small $\eps =0.025$}
		\put(52,0){\small $\eps =0.05$}
		\put(78,0){\small $\eps =0.075$}
		\put(2,62){\small $K_A$, $K_\eps/\rho_\eps$}
		\put(52,62){\small $C_A^{2323},C_\eps^{2323}/\rho_\eps$}
	\end{overpic}
%	\vspace{-3mm}
	\caption{
		The objective evolution of the derived shell lattices with different thickness $(2\eps)$ during the optimization, normalized by the volume fraction ($\rho_\eps$). `Asym.' represents the objective evaluated from  $\mathbf C_A$.
	}
	%    \vspace{-3mm}
	\label{fig:effect-evolve}
\end{figure}
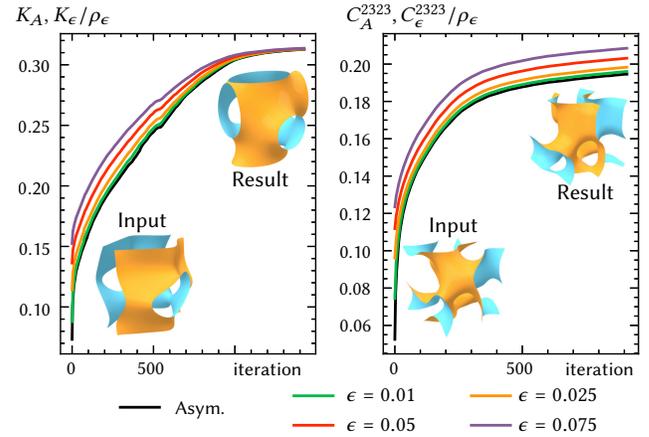

To validate whether optimizing ADS of the middle surfaces improves the effective stiffness of the shell lattice, we conduct experiments maximizing $K_A$ and $C^{2323}_A$.
The evolution of the objectives during optimization, for both the asymptotic and effective stiffness of the derived shell lattices with varying thicknesses, is shown in Figure~\ref{fig:effect-evolve}.
The results demonstrate that the effective stiffness exhibits the same trend as the asymptotic stiffness during the optimization.
Performance improvements for other optimized shell lattices in this work are summarized in Table~\ref{table:stat}.

\subsection{Comparisons}
A major distinction from traditional approaches based on shell elements is that our formulation and optimization of ADS focus exclusively on the geometry of the middle surface,  independent of the shell's thickness.
By definition, shell lattices optimized for higher ADS exhibit high, albeit not necessarily optimal, performance across a range of thicknesses with theoretical guarantees.
This  makes our approach particularly suitable for design scenarios where the thickness remains undetermined.
In contrast, the optimization that uses shell elements with a specific thickness introduces uncertainty regarding stiffness when the thickness varies, thus multiple optimizations are necessary.

Another advantage lies in our discretization scheme.
A typical  Reissner-Mindlin plate element used in non-parametric design~\cite{Shimoda2014-ao} assigns five DoFs to each node, while our method only assigns three DoFs.
This reduces the size of the linear system to be solved in the analysis step, thus decreasing computational costs.
A similar issue arises when resorting to  subdivision surfaces~\cite{BANDARA2016510,Cirak2010,BANDARA201862}, where numerical integration is much more expensive.
For parametric design methods~\cite{MA2022110426,Ramm01011993,Bletzinger1993-py} based on NURBS or B\'ezier patches, it becomes difficult to handle surfaces with complex topology, not only for the generation of the quadrilateral layout on a periodic surface, but also for addressing  the continuity problem  between patches and necking singularities that could appear during the optimization.

	\begin{remark}
		As thickness approaches zero, the higher-order infinitesimal bending energy in facet-based shell elements vanishes much more rapidly than the membrane strain energy,  making the elements degenerate into plane (membrane) elements. This essentially reduces the necessary  DoFs for evaluating strain energy. However, a simple plane element fails to capture surface deformation, as shown in Figure~\ref{fig:plate-singular} middle. Therefore, our scheme~\eqref{eq:gam-disc-final} could be considered as an improved membrane element, retaining less DoFs than shell elements while capable of characterizing  curved surfaces and their deformation.
	\end{remark}

\begin{table}[t]
	\centering
	\caption{
		Improvement of effective stiffness across different thicknesses ($2\eps$).
		$H^r_{2\eps}$ denotes the relative increase of the objective $f(\mathbf C_\eps/\rho_\eps)$, normalized by the corresponding increase in the objective $f(\mathbf C_A)$ after optimization. 
		For the surfaces in Figure~\ref{fig:cusobj}, the increase is computed without  the penalty term.
	} \label{table:stat}
%	\vspace{-3mm}
	\resizebox{0.99\linewidth}{!}
	{
		\begin{tabular}{lcccr|lcccr}
			\toprule
			Surface  & $H_{0.02}^r$ &$H_{0.05}^r$ & $H_{0.1}^r$ & $H_{0.15}^r$  &
			Surface  & $H_{0.02}^r$ &$H_{0.05}^r$ & $H_{0.1}^r$ & $H_{0.15}^r$  
			\\
			\midrule
			Fig.~\ref{fig:ads-max-G} [100] & 97\% & 93\% & 90\% & 87\%&
			Fig.~\ref{fig:abm-opt} (a)    & 96\% & 87\% & 78\% & 71\%  \\
			Fig.~\ref{fig:ads-max-G} [010] & 96\% & 89\% & 82\% & 77\% &
			Fig.~\ref{fig:abm-opt} (b)    & 89\% & 78\% & 69\% & 61\% \\
			Fig.~\ref{fig:ads-max-G} [001]  & 93\% & 84\% & 77\% & 72\% &
			Fig.~\ref{fig:abm-opt} (c)    & 93\% & 84\% & 75\% & 70\% \\
			Fig.~\ref{fig:ads-max-G} [011]  & 96\% & 91\% & 88\% & 85\% &
			Fig.~\ref{fig:abm-opt} (d)    & 90\% & 76\% & 64\% & 56\% \\
			Fig.~\ref{fig:ads-max-G} [101]  & 97\% & 93\% & 89\% & 87\%&
			Fig.~\ref{fig:vari-strain} (a)   & 96\% & 89\% & 82\% & 78\% \\
			Fig.~\ref{fig:ads-max-G} [110]  & 98\% & 96\% & 94\% & 92\%&
			Fig.~\ref{fig:vari-strain} (b)   & 96\% & 89\% & 83\% & 80\% \\
			Fig.~\ref{fig:ads-max-G} [111]  & 97\% & 95\% & 93\% & 91\%&
			Fig.~\ref{fig:vari-strain} (c)   & 89\% & 77\% & 68\% & 63\%  \\
			Fig.~\ref{fig:shear-opt} (a)  & 86\% & 70\% & 58\% & 52\% &
			Fig.~\ref{fig:vari-strain} (d)   & 93\% & 87\% & 82\% & 77\% \\
			Fig.~\ref{fig:shear-opt} (b)  & 89\% & 80\% & 73\% & 69\% &
			Fig.~\ref{fig:dev-strain} (a) & 87\% & 67\% & 52\% & 43\%\\
			Fig.~\ref{fig:shear-opt} (c)  & 95\% & 89\% & 85\% & 83\% &
			Fig.~\ref{fig:dev-strain} (b) & 73\% & 55\% & 48\% & 41\%\\
			Fig.~\ref{fig:shear-opt} (d)  & 92\% & 86\% & 80\% & 76\% &
			Fig.~\ref{fig:cusobj} (a) &77\% & 57\% & 42\% & 34\%\\
			Fig.~\ref{fig:shear-opt} (e)  & 99\% & 96\% & 93\% & 90\% &
			Fig.~\ref{fig:cusobj} (b) &77\% & 58\% & 44\% & 36\%\\
			Fig.~\ref{fig:shear-opt} (f)  & 97\% & 93\% & 89\% & 86\% &
			Fig.~\ref{fig:cusobj} (c) & 77\% & 58\% & 44\% & 36\%\\
			Fig.~\ref{fig:app-fence}      &90\% & 78\% & 70\% & 66\% &&&&&\\
			
			\bottomrule
		\end{tabular}
	}
\end{table}  

\section{Conclusion}
\label{sec:conc}
In this work, we introduce  asymptotic directional stiffness (ADS) to quantify the influence of the middle surfaces to the stiffness of shell lattices.
We establish the convergence theorem that provides its formulation.
Based on this formulation, we derive an  upper bound of ADS and establish the necessary and sufficient condition to achieve the bound.
This condition is applied to justify the optimal bulk modulus of TPMS shell lattice metamaterials.
To construct surfaces with higher ADS in a general case,  we propose numerical methods for evaluating and optimizing the ADS of 
periodic surfaces.
The simulation results validate our theoretical prediction and demonstrate the effectiveness of both the discretization and  optimization algorithms.

\paragraph{Limitation and future work}

Our optimization pipeline and implementation are not optimized for computational performance.
Acceleration techniques, such as Nesterov or Anderson acceleration~\cite{Aspremont_2021}, could be integrated into our algorithm.
However, as the topology of the mesh changes during the optimization, such application is not straightforward.
Additionally, more efficient solvers, such as geometric multigrid solvers~\cite{tut-multigrid}, could be used to solve the linear systems more efficiently.

The surgery operation is limited to simplifying the topology of the input surface.
As a result, the final optimized geometry is dependent on the initial topology; for example, our method cannot generate a high-genus surface from a low-genus initial shape. 
We anticipate that other topological operations beyond surgery could be explored in future work.
An alternative approach is using implicit representations for the surface, such as level set representation, which are more flexible to detect and handle topological variations.

Self-intersections have been observed in the experiments, which not only violates our assumption on the surface (embedding condition), but also introduces manufacturing problems, such as the creation of closed chambers.
In future work, we plan to incorporate repulsive energy terms as proposed in~\cite{Yu2021-srf,Sassen:2024:RS} to solve the problem.
Other manufacturing constraints like overhang angle constraint should also be taken into account.

As indicated by the formulation of ADS, it does not capture the bending stiffness of shell lattices.
	Pure bending deformations are essentially located in the inextensional displacement space ($\ker\gamma$).
	Shell structures are typically much less resistant to such deformations, as bending energy is of higher-order infinitesimal compared to membrane energy, and bending deformation induces significant stress near the boundary surfaces of shell.
	To identify these unstable deformation modes, one can solve the following eigenvalue problem
	\begin{equation}
		\min_{\boldsymbol{u}\in\boldsymbol H^1(\tilde\omega)\cap\mathbf{Rig}^\bot}\int_{\tilde\omega}\gamma(\boldsymbol{u}):\boldsymbol A_{\tilde\omega}:\gamma(\boldsymbol u)\quad{}s.t.\quad\|\boldsymbol u\|_{L^2(\tilde\omega)}^2=1
	\end{equation}
	based on our discrete scheme for \(\gamma(\boldsymbol{u})\).
	Consequently,  maximizing this minimum eigenvalue provides a possible strategy to reduce bending effects, and enhance buckling strength.

\begin{acks}
	We thank the anonymous reviewers for their constructive feedback. We are also grateful to Xiaoming Fu, Xiaoya Zhai, and Qing Fang for insightful discussions on potential applications of the theory, and to Tian Wu for reviewing the proof in an early draft.
	This  work  is  supported  by  the  National Key  R\&D  Program  of  China  (2022YFB3303400)  and  the  National Natural Science Foundation of China (62025207).
\end{acks}

\bibliographystyle{ACM-Reference-Format}
\bibliography{main}

\appendix

\section{Derivation of $E_M$}
\label{app:compact-form-EM}
It suffices to derive the following equation
\begin{equation}
	\label{eq:deri-Em-P}
	 E_M(\tw;\veps)=\veps:\overline{\mathbb{P}}:\veps =\int_{\tw}e_{\tw}:\A_{\tw}:e_{\tw}.
\end{equation}
The strain $e_{\tw}$ can be written in matrix form as:
\begin{equation}
	e_{\tw} = P \veps P^\top,
\end{equation}
where $P=\mathbf{I}-\n\n^\top$ is the projection matrix.
Then the integrand in~\eqref{eq:deri-Em-P}, denoted as $f_M(\tw;\veps)$ for brevity, becomes
\begin{equation}
	\begin{aligned}
		f_M(\tw;\veps)&=\lambda_0 P^{ik}P^{il}\veps_{kl}P^{jm}P^{jn}\veps_{mn}+2\mu P^{ik}P^{jl}\veps_{kl}P^{im}P^{jn}\veps_{mn}\\
		&=\veps:\mathbb P :\veps,
	\end{aligned}
\end{equation}
where
\begin{equation}
	\mathbb P:=(\lambda_0 P^{ik}P^{il}P^{jm}P^{jn}+2\mu P^{ik}P^{jl}P^{im}P^{jn})\e_k\otimes\e_l\otimes\e_m\otimes\e_n.
\end{equation}
Note that we are using the orthogonal basis, thus we sum over the duplicated indices, regardless of whether they are superscripts or subscripts.
Given that $P$ is a projection matrix, we have $P^\top P=P$, or $P^{ik}P^{jl}=P^{kl}$ in component form. Therefore, $\mathbb P$ is simplified as
\begin{equation}
	[\mathbb P]^{ijkl}=\lambda_0 P^{ij}P^{kl}+\mu\pr{P^{il}P^{jk}+P^{ik}P^{jl}}.
\end{equation}
Since $\veps$ is a constant tensor, we have
\begin{equation}
	E_M(\tw;\veps)=\frac{1}{|{\tw}|}\int_{\tw}\veps:\mathbb{ P}:\veps=\veps:\overline{\mathbb{P}}:\veps,
\end{equation}
where
$
\overline{\mathbb{P}}:=\frac{1}{|{\tw}|}\int_{\tw}\mathbb{ P}
$
denotes the average of $\mathbb P$ over $\tw$.

\section{Proof of Theorem~\ref{thm:eigen-sum-of-P}}
\label{app:eigen-sum-of-P}
\begin{proof}
	The sum of eigenvalues of $\mathbb P$ is
	\begin{equation}
		\label{eq:eigen-sum-P}
		\sum_{i=1}^6\lambda_i(\mathbb P)=\tr \mathbb P=\mathbb P :: \mathbb I,
	\end{equation}
	where $\mathbb I$ is the fourth-order identity tensor:
	\begin{equation}
		\mathbb I = \delta_{ik}\delta_{jl}\e^i\otimes\e^j\otimes\e^k\otimes\e^l.
	\end{equation}
	By computation, we have
	\begin{equation}
		{
			\mathbb P ::\mathbb I = \lambda_0 P_{kl}P_{kl} + \mu\pr{P_{kl}P_{kl}+P_{kk}P_{ll}}.
		}
	\end{equation}
	Notice that
	\begin{equation}
		\label{eq:trace-of-P}
		\begin{aligned}
			P_{kl}P_{kl}=\tr PP^\top =\tr P =2\\
			P_{kk}P_{ll}=\pr{\tr P}^2=4,
		\end{aligned}
	\end{equation}
	where the last equalities of both equations hold because $\tr P=\tr\pr{I-\n\n^\top}=2$.
	Hence, we get
	\begin{equation}
		\label{eq:p-double-dot-i}
		\mathbb P::\mathbb I = 2\lambda_0+6\mu.
	\end{equation}
	Then  we have
	\begin{equation}
		\sum_{i=1}^6 E_M(\tw;\veps_i)=\frac{1}{|\tw|}\int_{\tw}\sum_i\veps_i:\mathbb P:\veps_i=\frac{1}{|\tw|}\int_{\tw}\sum_{i=1}^6\lambda_i(\mathbb P)=2\lambda_0+6\mu,
	\end{equation}
	where we have utilized~\eqref{eq:eigen-sum-P} and~\eqref{eq:p-double-dot-i}.
\end{proof}

\section{Assemble the stiffness matrix and load vector}
\label{eq:assemb-stif-matvec}
We first select a 2D orthonormal frame 
\begin{equation}
	\label{eq:pf-defin}
	P_f=(\mathbf a_1,\mathbf a_2)^\top
\end{equation}
on each face as the coordinate basis to facilitate computation.
The strain and stress tensor are stored in the Voigt representation~\cite{voigt1910lehrbuch} under this basis:
\begin{equation}
	[e_f] = (e^{11}_f, e^{22}_f, 2 e^{12}_f)^\top.
\end{equation}
The second fundamental form $\bs b_f$ is approximated as the solution of the following equations for each edge vector $\bs l_{ij}$ of the face $f$: 
\begin{equation}
	\bs l_{ij}^\top \bs b_f\bs l_{ij} = -(\bs n_i-\bs n_j)\cdot\bs l_{ij}   .
\end{equation}
Once $\bs b_f$ is solved, the membrane strain is represented by 
\begin{equation}
	[\gamma(\bs u)]= (B_t - B_n)\bs u_f,
\end{equation}
where $\bs u_f=(\bs u_1;\bs u_2;\bs u_3)\in\mathbb R^{9}$ stacks the three vertex displacement vectors, 
and $B_t$ is the strain-displacement matrix mapping $\bs u_f$ to the Voigt representation of $\text{Sym}[\nabla\bs u_t]$ with $\bs u_t$ given by
\begin{equation}
	\bs u_t = \sum_i\phi_i P_f R_i^f\u_i .
\end{equation}
The matrix $B_n$ has the following form
\begin{equation}
	B_n = [\bs b_f] (\phi_1\bs n_1^\top,\phi_2\bs n_2^\top,\phi_3\bs n_3^\top),
\end{equation}
where $[\bs b_f]=(b_f^{11},b_f^{22},2b_f^{12})^\top$ and $\phi_i,\bs n_i$ denote the linear interpolation function (i.e., barycentric {coordinate}) and vertex normal, respectively.
Then the element stiffness matrix for each face is
\begin{equation}
	\label{eq:element-kf}
	\mathbf K_f = \int_{f} [\gamma(\bs u)]^\top D [\gamma(\bs u)],
\end{equation}
where $D$ is the membrane elastic matrix given by
\begin{equation}
	D=\begin{pmatrix}
		\lambda_0 + 2\mu &\lambda_0 & 0\\
		\lambda_0 &\lambda_0 + 2\mu & 0\\
		0&0 &\mu
	\end{pmatrix}.
\end{equation}
We use three point quadrature~\cite{XIAO2010663} on triangles to evaluate this integral. 
The load vector on each face is discretized as
\begin{equation}
	\label{eq:element-ff}
	\mathbf {f}_f =-\int_f [\gamma(\u)]^\top D [P_f \veps P_f^\top].
\end{equation}
Here, $[P_f \veps P_f^\top]$ denotes the Voigt representation of the projected strain $P_f\veps P_f^\top$.
This integral is evaluated using single point quadrature.
The element stiffness matrices in~\eqref{eq:element-kf} and the load vectors in~\eqref{eq:element-ff} are then assembled into the global stiffness matrix $\mathbf K^h$ and load vector $\mathbf f^h$.

\section{Discretization of sensitivity}
\label{app:sens}
We discretize the normal velocity on the mesh as a linear interpolation from the normal velocity on vertex, given by
\begin{equation}
	v_n^h :=\sum_{i\in\mathcal V}v_n^i\phi_i,
\end{equation}
where $v_n^i$ is the normal velocity on $i$-th vertex. 
$\phi_i$ is the linear interpolation function.

\paragraph{Sensitivity of Area}
The discretization of~\eqref{eq:area-rate} is given by
\begin{equation}
	\label{eq:area-time-rate}
	\dt{A}\approx-\sum_{f\in\mathcal F}\int_{f} v_n^h \tr{\bs b_f} = -\sum_{f\in\mathcal{F}}\sum_{i\in N(f)}v_n^i\int_{f}\phi_i\tr{\bs b_f},
\end{equation}
where $N(f)$ denotes the vertices of face $f$.
We use single point quadrature on each face to evaluate the integral.

\paragraph{Sensitivity of ADS}
The time derivative $\dt{I}_a^s$ is  discretized  as follows 
\begin{equation}
	\label{eq:Is-grad-h}
	\begin{aligned}
		\dt{I}^s_a\approx\sum_{f\in\mathcal{F}}&\sum_{i\in N(f)}v_n^i\int_{f} 2\phi_i\bigg[ B^h(\u^h,\u^h)-\frac{\tr \bs b_f}{2} A^h(\u^h,\u^h)
		\bigg]+ 2 C_i^h(\u^h,\u^h) .
	\end{aligned}
\end{equation}
Here, the symbols are computed as 
\begin{equation}
	\label{eq:dis-symb-def}
	\begin{aligned}
			A^h(\u^h,\v^h) &= e_f(\u^h):\A_{\tw}:e_f(\v^h)\\
			B^h(\u^h,\v^h)&=e_f(\u^h):\bs B_{\tw}:e_f(\v^h)\\
			C_i^h(\u^h,\v^h)&=e_f(\u^h):\A_{\tw}:\zeta_i^h(\v^h),
	\end{aligned}
\end{equation}
where  the strain $e_f(\u^h)$ is given by
\begin{equation}
	e_f(\u^h) = \text{Sym}[\u_t^h]- u_3^h\bs b_f+P_f\veps P_f^\top,
\end{equation}
and $P_f$ denotes the projection given in~\eqref{eq:pf-defin}. The term $\zeta^h_i(\bar\u^h)$ is defined as
\begin{equation}
	\begin{aligned}
		\zeta_i^h(\u^h)&=\bs b_f\u_t^h (\nabla \phi_i)^\top-\phi_i\bs b_f \nabla\u_t^h+\nabla u_3^h(\nabla \phi_i)^\top + u_3^h \phi_i\bs{c}_f\\&\quad+2P_f\veps\n_f(\nabla \phi_i)^\top,
	\end{aligned}
\end{equation}
where $\n_f$ denotes the normal vector of face $f$ and $\bs c_f=\bs b_f\cdot\bs b_f$.
The tensor contraction in~\eqref{eq:dis-symb-def} is evaluated in matrix arithmetic:
\begin{equation}
	\begin{aligned}
			e_f^1:\A_{\tw}:e_f^2 &= \lambda_0 \tr{e_f^1}\tr{e_f^2}+2\mu \tr(e_f^1 e_f^2)\\
			e_f^1:\bs B_{\tw}:e_f^2&=\lambda_0 \pr{\tr(\bs{b}_f e_f^1)\tr e_f^2+\tr(\bs{b}_f e_f^2)\tr e_f^1} + 4\mu\tr(\bs b_f e_f^1 e_f^2),
	\end{aligned}
\end{equation}
where $e_f^1, e_f^2$ denote arbitrary symmetric matrices on face $f$.

Similarly, the discretization of $\dt{I}^s_{ijkl}$ in~\eqref{eq:ca-ijkl-dt} is given by 
\begin{equation}
	\begin{aligned}
		\dt{I}^s_{ijkl}\approx\sum_{f\in\mathcal{F}}&\sum_{q\in N(f)}v_n^q\int_{f} 2\phi_q\bigg[ B^h(\u_{ij}^h,\u_{kl}^h)-\frac{\tr \bs b_f}{2} A^h(\u_{ij}^h,\u_{kl}^h)
		\bigg]+ \\ &C_q^h(\u_{ij}^h,\u_{kl}^h)  + C_q^h(\u^h_{kl},\u^h_{ij}).
	\end{aligned}
\end{equation}

\section{Isotropic penalty}
\label{sec:iso-penal}
For cubic symmetric shell lattice, the isotropic condition is typically satisfied by minimizing or constraining the error between Zener ratio~\cite{zener1948elasticity} and 1.
Since the input surface is generally not cubic symmetric, we adopt a different approach.
Note that the isotropic elastic tensors have the following form
\begin{equation}
	\mathbb C = \lambda \delta^{ij}\delta^{kl}+\mu\pr{\delta^{ik}\delta^{jl}+\delta^{il}\delta^{jk}}.
\end{equation}
Therefore, they constitute a two-dimensional linear space spanned by two tensors
\begin{equation}
	\mathbb C_\lambda := \delta^{ij}\delta^{kl},\quad \mathbb C_\mu := \delta^{ik}\delta^{jl}+\delta^{il}\delta^{jk}.
\end{equation}
We penalize the anisotropy of ADS by adding the following penalty term to the objective 
\begin{equation}
	\label{eq:els-iso-pen}
	f^s_{iso}(\mathbb C_A) = \min_{\lambda,\mu}\, \frac{1}{2}\|\mathbb C_A-(\lambda\mathbb C_\lambda +\mu\mathbb C_\mu) \|_F^2.
\end{equation}
The time derivative of this function is given by
\begin{equation}
	\frac{d}{dt}f^s_{iso}(\mathbb C_A) = [(\lambda^*\mathbb C_\lambda+\mu^*\mathbb C_\mu)  -\mathbb C_A]:: \dt{\mathbb C}_A ,
\end{equation}
where $\lambda^*,\mu^*$ denote the minimizer of~\eqref{eq:els-iso-pen} and `::' denotes tensor contraction.

\section{Objective for negative Poisson's ratio}
\label{sec:opt-npr}
We minimize the following objective function to design shell lattices with negative Poisson's ratios
\begin{equation}
	f(\mathbf C_A) = \eps_1(\mathbf C_A)+\eps_2(\mathbf C_A) .
\end{equation}
Here, $\eps_1, \eps_2$ denote the lateral contractions in two orthogonal  directions when the material is stretched. We approximate them by solving the following problem based on the the principle of minimal potential energy:
\begin{equation}
	\label{eq:ps-def-prob}
	\min_{\eps_1,\eps_2}\,F(\eps_1,\eps_2):\mathbf C_A: F(\eps_1,\eps_2),
\end{equation}
where
\begin{equation}
	F(\eps_1,\eps_2) =
	\begin{pmatrix}
		-\eps_1 & 0 & 0 \\
		0        & -\eps_2 & 0 \\
		0        & 0 & 1 \\
	\end{pmatrix}.
\end{equation}
Equation~\eqref{eq:ps-def-prob} is a least square problem, whose solution can be readily found by solving a linear system.

\end{document}